\documentclass[english,journal]{IEEEtran}
\usepackage[T1]{fontenc}
\usepackage{xcolor}
\usepackage{float}
\usepackage{amsbsy}
\usepackage{amstext}
\usepackage{amsthm}
\usepackage{graphicx}

\makeatletter

\floatstyle{ruled}
\newfloat{algorithm}{tbp}{loa}
\providecommand{\algorithmname}{Algorithm}
\floatname{algorithm}{\protect\algorithmname}

\theoremstyle{plain}
\newtheorem{thm}{\protect\theoremname}
\theoremstyle{remark}
\newtheorem{rem}[thm]{\protect\remarkname}

\usepackage{amsmath}
\usepackage{amssymb}
\usepackage[numbers,sort]{natbib}

\usepackage{graphicx}
\usepackage{subfigure}
\usepackage{pdfsync}
\usepackage{url}
\usepackage{pdfsync}
\usepackage{todonotes}
\makeatother

\usepackage{babel}
\providecommand{\remarkname}{Remark}
\providecommand{\theoremname}{Theorem}

\begin{document}
\title{Inexact Block Coordinate Descent Algorithms for Nonsmooth Nonconvex
Optimization}
\author{Yang Yang, Marius Pesavento, Zhi-Quan Luo, Björn Ottersten\thanks{Y. Yang is with the Competence Center for High Performance Computing,
Fraunhofer Institute for Industrial Mathematics, 67663 Kaiserslautern,
Germany (email: yang.yang@itwm.fraunhofer.de). His work is supported
by the ERC project AGNOSTIC.}\thanks{M. Pesavento is with the Communication Systems Group, Technische Universität
Darmstadt, 64283 Darmstadt, Germany (email: pesavento@nt.tu-darmstadt.de).
His work is supported by supported by the EXPRESS project within the
DFG priority program CoSIP (DFG-SPP 1798).}\thanks{Z.-Q. Luo is with Shenzhen Research Institute of Big Data, and the
Chinese University of Hong Kong, Shenzhen, China (email: luozq@cuhk.edu.cn).
His work is supported by the leading talents of Guangdong province
Program (No. 00201501), the National Natural Science Foundation of
China (No. 61731018), the Development and Reform Commission of Shenzhen
Municipality, and the Shenzhen Fundamental Research Fund (No. KQTD201503311441545).}\thanks{B. Ottersten is with University of Luxembourg, L-1855 Luxembourg (email:
bjorn.ottersten@uni.lu). His work is supported by the ERC project
AGNOSTIC.}}
\maketitle
\begin{abstract}
In this paper, we propose an inexact block coordinate descent algorithm
for large-scale nonsmooth nonconvex optimization problems. At each
iteration, a particular block variable is selected and updated by
inexactly solving the original optimization problem with respect to
that block variable. More precisely, a local approximation of the
original optimization problem is solved. The proposed algorithm has
several attractive features, namely, i) high flexibility, as the approximation
function only needs to be strictly convex and it does not have to
be a global upper bound of the original function; ii) fast convergence,
as the approximation function can be designed to exploit the problem
structure at hand and the stepsize is calculated by the line search;
iii) low complexity, as the approximation subproblems are much easier
to solve and the line search scheme is carried out over a properly
constructed differentiable function; iv) guaranteed convergence of
a subsequence to a stationary point, even when the objective function
does not have a Lipschitz continuous gradient. Interestingly, when
the approximation subproblem is solved by a descent algorithm, convergence
of a subsequence to a stationary point is still guaranteed even if
the approximation subproblem is solved inexactly by terminating the
descent algorithm after a finite number of iterations. These features
make the proposed algorithm suitable for large-scale problems where
the dimension exceeds the memory and/or the processing capability
of the existing hardware. These features are also illustrated by several
applications in signal processing and machine learning, for instance,
network anomaly detection and phase retrieval.
\end{abstract}

\begin{IEEEkeywords}
Big Data, Block Coordinate Descent, Phase Retrieval, Line Search,
Network Anomaly Detection, Successive Convex Approximation
\end{IEEEkeywords}

\section{Introduction}

In this paper, we consider the optimization problem
\begin{align}
\underset{\mathbf{x}=(\mathbf{x}_{k})_{k=1}^{K}}{\textrm{minimize}}\quad & h(\mathbf{x})\triangleq f(\mathbf{x}_{1},\ldots,\mathbf{x}_{K})+\underbrace{\sum_{k=1}^{K}g_{k}(\mathbf{x}_{k})}_{g(\mathbf{x})},\nonumber \\
\textrm{subject to}\quad & \mathbf{x}_{k}\in\mathcal{X}_{k}\subseteq\mathbb{R}^{I_{k}},\forall k=1,\ldots,K,\label{eq:problem-formulation}
\end{align}
where the function $h$ is proper, $f$ is smooth (but not necessarily
convex), $g_{k}$ is proper, lower semicontinuous and convex (but
not necessarily smooth), and the constraint set has a Cartesian product
structure with $\mathcal{X}_{k}$ being closed and convex for all
$k=1,\ldots,K$. Such a formulation plays a fundamental role in signal
processing and machine learning, and typically $f$ models the estimation
error or empirical loss while $g_{k}$ is a regularization (penalty)
function promoting in the solution a certain structure known a priori
such as sparsity.

For a large-scale nonconvex optimization problem of the form (\ref{eq:problem-formulation}),
the block coordinate descent (BCD) algorithm has been recognized as
an efficient and reliable numerical method. Its variable update is
based on the so-called nonlinear best-response \cite{Bertsekas,Tseng2001,Razaviyayn2013,Beck2013,Wright2015}:
at each iteration of the BCD algorithm, one block variable, say $\mathbf{x}_{k}$,
is updated by its best-response while the other block variables are
fixed to their values of the preceding iteration
\begin{align}
\mathbf{x}_{k}^{t+1} & =\underset{\mathbf{x}_{k}\in\mathcal{X}_{k}}{\arg\min}\;h(\mathbf{x}_{1}^{t+1},\ldots,\mathbf{x}_{k-1}^{t+1},\mathbf{x}_{k},\mathbf{x}_{k+1}^{t},\ldots,\mathbf{x}_{K}^{t})\label{eq:BCD}\\
 & =\underset{\mathbf{x}_{k}\in\mathcal{X}_{k}}{\arg\min}\;f((\mathbf{x}_{j}^{t+1})_{j=1}^{k-1},\mathbf{x}_{k},(\mathbf{x}_{j}^{t})_{j=k+1}^{K})+g_{k}(\mathbf{x}_{k}).\nonumber
\end{align}
That is, the best-response is the optimal point that minimizes $h(\mathbf{x})$
with respect to (w.r.t.) the variable $\mathbf{x}_{k}$.

The BCD algorithm has several notable advantages. First of all, the
subproblem (\ref{eq:BCD}) (w.r.t. a block variable $\mathbf{x}_{k}$)
is much easier to solve than the original problem (\ref{eq:problem-formulation})
(w.r.t. the whole set of variables $\mathbf{x}$), and the best-response
even has a closed-form expression in many applications, for example
LASSO \cite{Bach2012}. It is thus suitable for implementation on
hardware with limited memory and/or computational capability. Secondly,
as all block variables are updated sequentially, when a block variable
is updated, the newest value of other block variables is always incorporated.
These two attractive features can sometimes lead to even faster convergence
than their parallel counterpart, namely, the Jacobi algorithm (also
known as the parallel best-response algorithm) \cite{Bertsekas}.

In cases where the subproblems (\ref{eq:BCD}) are still difficult
to solve and/or (sufficient) convergence conditions (mostly on the
convexity of $f$ and the uniqueness of $\mathbf{x}_{k}^{t+1}$, see
\cite{bertsekas1999nonlinear,Tseng2001,Xu2013a} and the references
therein) are not satisfied, several extensions have been proposed.
Their central idea is to solve the optimization problem (\ref{eq:BCD})
\emph{inexactly}. For example, in the block successive upper bound
minimization (BSUM) algorithm \cite{Razaviyayn2013}, a global upper
bound function of $h((\mathbf{x}_{j}^{t+1})_{j=1}^{k-1},\mathbf{x}_{k},(\mathbf{x}_{j}^{t})_{j=k+1}^{K})$
is minimized at each iteration. Common examples are proximal approximations
\cite{Xu2013a,Chouzenoux2016} and, if $\nabla f$ is block Lipschitz
continuous\footnote{The gradient $\nabla f$ is Lipschitz continuous if there exists a
finite constant $L$ such that $\left\Vert \nabla f(\mathbf{x})-\nabla f(\mathbf{y})\right\Vert \le L\left\Vert \mathbf{x}-\mathbf{y}\right\Vert $
for all $\mathbf{x},\mathbf{y}\in\mathcal{X}_{1}\times\ldots\times\mathcal{X}_{K}$.
It is block Lipschitz continuous if there exists a finite constant
$L$ such that $\left\Vert \nabla_{k}f(\mathbf{x}_{k},\mathbf{x}_{-k})-\nabla_{k}f(\mathbf{y}_{k},\mathbf{x}_{-k})\right\Vert \leq L\left\Vert \mathbf{x}_{k}-\mathbf{y}_{k}\right\Vert $
for all $\mathbf{x}_{k},\mathbf{y}_{k}\in\mathcal{X}_{k}$ and $\mathbf{x}_{-k}\in\mathcal{X}_{1}\times\ldots\mathcal{X}_{k-1}\times\mathcal{X}_{k+1}\times\ldots\times\mathcal{X}_{K}$
and $k=1,\ldots,K$.}, proximal-linear approximation \cite{Xu2013a,Bolte2014}. However,
for the BSUM algorithm, a global upper bound function may not exist
for some $f$ (and thus $h$).

The block Lipschitz continuity assumption is not needed if a stepsize
is employed in the variable update. In practice, the stepsize can
be determined by line search \cite{Mine1981,Tseng2009}. Nevertheless,
only a specific approximation of $f$ is considered, namely, quadratic
approximation. Sometimes it may be desirable to use other approximations
to better exploit the problem structure, for example, best-response
approximation and partial linearization approximation when the nonconvex
function $f$ has ``partial'' convexity (their precise descriptions
are provided in Section \ref{sec:Block-SCA}) . This is the central
idea in recent (parallel) successive convex approximation (SCA) algorithms
\cite{Patriksson1995,Scutarib,Scutari_BigData,Yang_ConvexApprox,Yang_NonconvexRegularization}
and block successive convex approximation (BSCA) algorithms \cite{Razaviyayn2013,Razaviyayn2014,Scutari_hybrid},
which consist in solving a sequence of successively refined convex
approximation subproblems. A new line search scheme to determine the
stepsize is also proposed in \cite{Yang_ConvexApprox,Yang_NonconvexRegularization}:
it is carried out over a properly constructed smooth function and
its complexity is much lower than traditional schemes that directly
operate on the original nonsmooth function \cite{Mine1981,Patriksson1995,Tseng2009}.
For example, as we will see later in the applications studied in this
paper, when $f$ represents a quadratic loss function, the exact line
search has a simple analytical expression.

Nevertheless, existing BSCA schemes also have their limitations: the
BSCA algorithm proposed in \cite{Razaviyayn2013} is not applicable
when the objective function is nonsmooth, and the convergence of the
BSCA algorithms proposed in \cite{Razaviyayn2014,Scutari_hybrid}
is only established under the assumption that $\nabla f$ is Lipschitz
continuous and the stepsizes are decreasing. Although it is shown
in \cite{Razaviyayn2014} that constant stepsizes can also be used,
the choice of the constant stepsizes depends on the Lipschitz constant
of $\nabla f$ that is not easy to obtain/estimate when the problem
dimension is extremely large.

The standard SCA and BSCA algorithms \cite{Tseng2009,Razaviyayn2013,Razaviyayn2014,Yang_ConvexApprox,Yang_NonconvexRegularization}
are based on the assumption that the approximation subproblem is solved
perfectly at each iteration. Unless the approximation subproblems
have a closed-form solution, this assumption can hardly be satisfied
by iterative algorithms that exhibit an asymptotic convergence only
as they must be terminated after a finite number of iterations in
practice. It is shown in \cite{Scutari_BigData,Scutari_hybrid,Chouzenoux2016}
that convergence is still guaranteed if the approximation subproblems
are solved approximately with a prescribed accuracy. However, the
solution accuracy is specified by an error bound which is difficult
to verify in practice. A different approach is adopted in \cite{Bonettini2011,Bonettini2018}
where the optimization problem (\ref{eq:BCD}) is solved inexactly
by running the (proximal) gradient projection algorithm for a finite
number of iterations. Nevertheless, its convergence is only established
for the specific application in nonnegative matrix factorization in
\cite{Bonettini2011} and the use of the (proximal) gradient projection
can be restrictive.

In this paper, we propose a block successive convex approximation
(BSCA) framework for the nonsmooth nonconvex problem (\ref{eq:problem-formulation})
by extending the parallel update scheme in \cite{Yang_ConvexApprox}
to a block update scheme. The proposed BSCA algorithm consists in
optimizing a sequence of successively refined approximation subproblems,
and has several attractive features.
\begin{enumerate}
\item [i)] The approximation function is a strictly convex approximation
of the original function and it does not need to be a global upper
bound of the original function;
\item [ii)] The stepsize is calculated by performing the (exact or successive)
line search scheme along the coordinate of the block variable being
updated and has low complexity as it is carried out over a properly
constructed smooth function;
\item [iii)] If the approximation subproblem does not admit a closed-form
solution and is solved iteratively by a descent algorithm, for example
the (parallel) SCA algorithm proposed in \cite{Yang_ConvexApprox},
the descent algorithm can be terminated after a finite number of iterations;
\item [iv)] Convergence of a subsequence to a stationary point is established,
even when $f$ is not multiconvex and/or $\nabla f$ is not block
Lipschitz continuous.
\end{enumerate}
These features are distinctive from existing works from the following
aspects:
\begin{itemize}
\item Feature i) extends the BSUM algorithm \cite{Razaviyayn2013} and BCD
algorithm \cite[(1.3b)]{Xu2013a} where the approximation function
must be a global upper bound of the original function, \cite{Mine1981,Tseng2009}
and \cite{Razaviyayn2014,Scutari_hybrid} where the approximation
functions must be quadratic and strongly convex, respectively;
\item Feature ii) extends \cite{Razaviyayn2014,Scutari_hybrid} where decreasing
stepsizes are used, and \cite{Mine1981,Patriksson1995,Tseng2009}
where the line search is over the original nonsmooth function and
has a high complexity;
\item Feature iii) extends \cite{Scutari_hybrid} where the approximation
subproblems must be solved with increasing accuracy. We remark that
this feature is inspired by \cite{Patriksson1995}, but we establish
convergence under weaker assumptions;
\item Feature iv) extends \cite[(1.3a)]{Xu2013a} where $f$ is multi-strongly-convex,
\cite{Razaviyayn2014,Scutari_hybrid,Chouzenoux2016} where $\nabla f$
must be Lipschitz continuous, \cite[(1.3c)]{Xu2013a} and \cite{Bolte2014}
where $\nabla f$ must be block Lipschitz continuous, and \cite{Patriksson1995,Tseng2009}
where line search over the original nonsmooth function is used. Nevertheless,
the convergence of a subsequence is weaker than the convergence of
the whole sequence established in \cite{Xu2013a,Chouzenoux2016,Bolte2014}.
\end{itemize}
These attractive features are illustrated by several applications
in signal processing and machine learning, namely, network anomaly
detection and phase retrieval.

The rest of the paper is structured as follows. In Sec. \ref{sec:Review-of-SCA},
we give a brief review of the SCA framework proposed in \cite{Yang_ConvexApprox}.
In Sec. \ref{sec:Block-SCA}, the BSCA framework together with the
convergence analysis is formally presented. An inexact BSCA framework
is proposed in Sec. \ref{sec:Inexact-Block-SCA}. The attractive features
of the proposed (exact and inexact) BSCA framework are illustrated
through several applications in Sec. \ref{sec:Example-applications}.
Finally some concluding remarks are drawn in Sec. \ref{sec:Concluding-Remarks}.

\emph{Notation: }We use $x$, $\mathbf{x}$ and $\mathbf{X}$ to denote
a scalar, vector and matrix, respectively. We use $x_{j,k}$ and $\mathbf{x}_{j}$
to denote the $(j,k)$-th element and the $j$-th column of $\mathbf{X}$,
respectively; $x_{k}$ is the $k$-th element of $\mathbf{x}$ where
$\mathbf{x}=(x_{k})_{k=1}^{K}$, and $\mathbf{x}_{-k}$ denotes all
elements of $\mathbf{x}$ except $x_{k}$: $\mathbf{x}_{-k}=(x_{j})_{j=1,j\neq k}^{K}$.
We denote $\mathbf{x}^{p}$ and $\mathbf{x/y}$ as the element-wise
operation, i.e., $(\mathbf{x}^{p})_{k}=(x_{k})^{p}$ and $(\mathbf{x}/\mathbf{y})_{k}=x_{k}/y_{k}$,
respectively. Notation $\mathbf{x}\circ\mathbf{y}$ denotes the Hadamard
product between $\mathbf{x}$ and $\mathbf{y}$. The operator $[\mathbf{x}]_{\mathbf{a}}^{\mathbf{b}}$
returns the element-wise projection of $\mathbf{x}$ onto $[\mathbf{a,b}]$:
$[\mathbf{x}]_{\mathbf{a}}^{\mathbf{b}}\triangleq\max(\min(\mathbf{x},\mathbf{b}),\mathbf{a})$.
We denote $\mathbf{d}(\mathbf{X})$ as the vector that consists of
the diagonal elements of $\mathbf{X}$ and $\textrm{diag}(\mathbf{x})$
is a diagonal matrix whose diagonal vector is $\mathbf{x}$. We use
$\mathbf{1}$ to denote a vector with all elements equal to 1. The
operator $\left\Vert \mathbf{X}\right\Vert _{p}$ specifies the $p$-norm
of $\mathbf{X}$ and it denotes the spectral norm when $p$ is not
specified. $S_{\mathbf{a}}(\mathbf{b})$ denotes the soft-thresholding
operator: $S_{\mathbf{a}}(\mathbf{b})\triangleq\max(\mathbf{b}-\mathbf{a},\mathbf{0})-\max(\mathbf{-b-a},0)$.

\section{\label{sec:Review-of-SCA}Review of the Successive Convex Approximation
framework}

In this section, we present a brief review of (a special case of)
the SCA framework developed in \cite{Yang_ConvexApprox} for problem
(\ref{eq:problem-formulation}). It consists of solving a sequence
of successively refined approximation subproblems: given $\mathbf{x}^{t}$
at iteration $t$, the approximation function of $f(\mathbf{x})$
w.r.t. $\mathbf{x}_{k}$ is denoted as $\widetilde{f}_{k}(\mathbf{x}_{k};\mathbf{x}^{t})$,
and the approximation subproblem consists of minimizing the approximation
function $\widetilde{h}(\mathbf{x};\mathbf{x}^{t})\triangleq\sum_{k=1}^{K}\widetilde{f}_{k}(\mathbf{x}_{k};\mathbf{x}^{t})+\sum_{k=1}^{K}g_{k}(\mathbf{x}_{k})$
over the constraint set $\mathcal{X}_{1}\times\ldots\times\mathcal{X}_{K}$:
\begin{equation}
\mathbb{B}\mathbf{x}^{t}\in\underset{(\mathbf{x}_{k}\in\mathcal{X}_{k})_{k=1}^{K}}{\arg\min}\biggl\{\underbrace{\sum_{k=1}^{K}\widetilde{f}_{k}(\mathbf{x}_{k};\mathbf{x}^{t})}_{\widetilde{f}(\mathbf{x};\mathbf{x}^{t})}+\underbrace{\sum_{k=1}^{K}g_{k}(\mathbf{x}_{k})}_{g(\mathbf{x})}\biggr\},\label{eq:review-approximate-problem}
\end{equation}
where $\widetilde{f}_{k}(\mathbf{x}_{k};\mathbf{x}^{t})$ satisfies
several technical assumptions, most notably,
\begin{itemize}
\item Convexity\footnote{Please refer to \cite[Sec. II]{Yang_ConvexApprox} for optimization
terminologies such as (strict, strong) convexity, descent direction
and stationary point.}: The function $\widetilde{f}_{k}(\mathbf{x}_{k};\mathbf{x}^{t})$
is convex in $\mathbf{x}_{k}$ for any given $\mathbf{x}^{t}\in\mathcal{X}$;
\item Gradient Consistency: The gradient of $\widetilde{f}_{k}(\mathbf{x}_{k};\mathbf{x}^{t})$
and the gradient of $f(\mathbf{x})$ are identical at $\mathbf{x}=\mathbf{x}^{t}$,
i.e., $\nabla_{\mathbf{x}_{k}}\widetilde{f}_{k}(\mathbf{x}_{k}^{t};\mathbf{x}^{t})=\nabla_{\mathbf{x}_{k}}f(\mathbf{x}^{t})$.
\end{itemize}
We have also implicitly assumed that $\mathbb{B}\mathbf{x}^{t}$ exists.
The approximation subproblem (\ref{eq:review-approximate-problem})
can readily be decomposed into $K$ independent subproblems that can
be solved in parallel: $\mathbb{B}\mathbf{x}^{t}=(\mathbb{B}_{k}\mathbf{x}^{t})_{k=1}^{K}$
and
\[
\underset{\mathbf{x}_{k}\in\mathcal{X}_{k}}{\min}\bigl\{\widetilde{f}_{k}(\mathbf{x}_{k};\mathbf{x}^{t})+g_{k}(\mathbf{x}_{k})\bigr\},k=1,\ldots,K.
\]

\begin{rem}
The approximation function $\sum_{k=1}^{K}\widetilde{f}_{k}(\mathbf{x}_{k};\mathbf{x}^{t})$
in (\ref{eq:review-approximate-problem}) is a special case of the
general SCA framework developed in \cite{Yang_ConvexApprox} because
it is separable among the different block variables. More generally,
$\widetilde{f}(\mathbf{x};\mathbf{x}^{t})$, the approximation function
of $f(\mathbf{x})$, only needs to be convex and differentiable with
the same gradient as $f(\mathbf{x})$ at $\mathbf{x}^{t}$, and it
does not necessarily admit a separable structure.
\end{rem}
Since $\mathbb{B}\mathbf{x}^{t}$ is an optimal point of problem (\ref{eq:review-approximate-problem}),
we have
\begin{align}
0 & \overset{(a)}{\geq}\widetilde{f}(\mathbb{B}\mathbf{x}^{t};\mathbf{x}^{t})+g(\mathbb{B}\mathbf{x}^{t})-(\widetilde{f}(\mathbf{x}^{t};\mathbf{x}^{t})+g(\mathbf{x}^{t}))\nonumber \\
 & \overset{(b)}{\geq}(\mathbb{B}\mathbf{x}^{t}-\mathbf{x}^{t})^{T}\nabla\widetilde{f}(\mathbf{x}^{t};\mathbf{x}^{t})+g(\mathbb{B}\mathbf{x}^{t})-g(\mathbf{x}^{t})\nonumber \\
 & \overset{(c)}{=}(\mathbb{B}\mathbf{x}^{t}-\mathbf{x}^{t})^{T}\nabla f(\mathbf{x}^{t})+g(\mathbb{B}\mathbf{x}^{t})-g(\mathbf{x}^{t})\triangleq d(\mathbf{x}^{t}),\label{eq:review-descent-direction}
\end{align}
where $(a)$, $(b)$ and $(c)$ is due to the optimality of $\mathbb{B}\mathbf{x}^{t}$,
the convexity of $\widetilde{f}(\mathbf{x};\mathbf{x}^{t})$ in $\mathbf{x}$
and the gradient consistency assumption, respectively. Therefore $\mathbb{B}\mathbf{x}^{t}-\mathbf{x}^{t}$
is a descent direction of the original objective function $h(\mathbf{x})$
in (\ref{eq:problem-formulation}) along which the function value
$h(\mathbf{x})$ can be further decreased compared with $h(\mathbf{x}^{t})$
\cite[Prop. 1]{Yang_ConvexApprox}. This motivates us to refine $\mathbf{x}^{t}$
and define $\mathbf{x}^{t+1}$ as follows:
\begin{equation}
\mathbf{x}^{t+1}=\mathbf{x}^{t}+\gamma^{t}(\mathbb{B}\mathbf{x}^{t}-\mathbf{x}^{t}),\label{eq:review-variable-update}
\end{equation}
where $\gamma^{t}\in(0,1]$ is the stepsize that needs to be selected
properly to yield a fast convergence.

It is natural to select a stepsize such that the function $h(\mathbf{x}^{t}+\gamma(\mathbb{B}\mathbf{x}^{t}-\mathbf{x}^{t}))$
is minimized w.r.t. $\gamma$:
\begin{equation}
\min_{0\leq\gamma\leq1}f(\mathbf{x}^{t}+\gamma(\mathbb{B}\mathbf{x}^{t}-\mathbf{x}^{t}))+g(\mathbf{x}^{t}+\gamma(\mathbb{B}\mathbf{x}^{t}-\mathbf{x}^{t})),\label{eq:review-traditional-line-search}
\end{equation}
and this is the so-called exact line search (also known as the minimization
rule). For nonsmooth optimization problems, the traditional exact
line search usually suffers from a high complexity as the optimization
problem (\ref{eq:review-traditional-line-search}) is nondifferentiable.
It is shown in \cite[Sec. III-A]{Yang_ConvexApprox} that the stepsize
obtained by performing the exact line search over the following \emph{differentiable}
function also yields a decrease in $h(\mathbf{x})$:
\begin{equation}
\gamma^{t}\in\underset{0\leq\gamma\leq1}{\textrm{argmin}}\bigl\{ f(\mathbf{x}^{t}+\gamma(\mathbb{B}\mathbf{x}^{t}-\mathbf{x}^{t}))+g(\mathbf{x}^{t})+\gamma(g(\mathbb{B}\mathbf{x}^{t})-g(\mathbf{x}^{t}))\bigr\}.\label{eq:review-proposed-stepsize-exact}
\end{equation}
To see this, we remark that firstly, the objective function in (\ref{eq:review-proposed-stepsize-exact})
is an upper bound of the objective function in (\ref{eq:review-traditional-line-search})
which is tight at $\gamma=0$ since $g$ is convex:
\[
g(\mathbf{x}^{t}+\gamma(\mathbb{B}\mathbf{x}^{t}-\mathbf{x}^{t}))\leq(1-\gamma)g(\mathbf{x}^{t})+\gamma g(\mathbb{B}\mathbf{x}^{t}),0\leq\gamma\leq1.
\]
Secondly, the objective function in (\ref{eq:review-proposed-stepsize-exact})
has a negative slope at $\gamma=0$ as its gradient is equal to $d(\mathbf{x}^{t})$
in (\ref{eq:review-descent-direction}). Therefore, $\gamma^{t}>0$
and $h(\mathbf{x}^{t}+\gamma^{t}(\mathbb{B}\mathbf{x}^{t}-\mathbf{x}^{t}))<h(\mathbf{x}^{t})$.

\begin{algorithm}[t]
\textbf{Initialization: }$t=0$ and $\mathbf{x}^{0}\in\mathcal{X}$
(arbitrary but fixed).

Repeat the following steps until convergence:

\begin{enumerate}

\item[\textbf{S1: }]Compute $\mathbb{B}\mathbf{x}^{t}=(\mathbb{B}_{k}\mathbf{x}^{t})_{k=1}^{K}$
by solving the following independent optimization problems in parallel:
\[
\mathbb{B}_{k}\mathbf{x}^{t}=\underset{\mathbf{x}_{k}\in\mathcal{X}_{k}}{\arg\min}\bigl\{\widetilde{f}_{k}(\mathbf{x}_{k};\mathbf{x}^{t})+g(\mathbf{x}_{k})\bigr\},\,k=1,\ldots,K.
\]

\item[\textbf{S2: }]Compute $\gamma^{t}$ by the exact line search
(\ref{eq:review-proposed-stepsize-exact}) or the successive line
search (\ref{eq:review-proposed-stepsize-successive}).

\item[\textbf{S3: }]Update $\mathbf{x}$: $\mathbf{x}^{t+1}=\mathbf{x}^{t}+\gamma^{t}(\mathbb{B}\mathbf{x}^{t}-\mathbf{x}^{t})$.

\item[\textbf{S4: }]$t\leftarrow t+1$ and go to \textbf{S1}.

\end{enumerate}

\caption{\label{alg:Successive-approximation-method}The parallel successive
convex approximation algorithm for nonsmooth nonconvex optimization
problem (\ref{eq:problem-formulation}) (proposed in \cite{Yang_ConvexApprox})}
\end{algorithm}

If the scalar differentiable optimization problem in (\ref{eq:review-proposed-stepsize-exact})
is still difficult to solve, the low-complexity successive line search
(also known as the Armijo rule) can be used instead \cite[Sec. III-A]{Yang_ConvexApprox}:
given scalars $0<\alpha<1$ and $0<\beta<1$, the stepsize $\gamma^{t}$
is set to be $\gamma^{t}=\beta^{m_{t}}$, where $m_{t}$ is the smallest
nonnegative integer $m$ satisfying
\begin{align}
 & f(\mathbf{x}^{t}+\beta^{m}(\mathbb{B}\mathbf{x}^{t}-\mathbf{x}^{t}))+g(\mathbf{x}^{t})+\beta^{m}(g(\mathbb{B}\mathbf{x}^{t})-g(\mathbf{x}^{t}))\nonumber \\
 & \leq f(\mathbf{x}^{t})+g(\mathbf{x}^{t})+\alpha\beta^{m}d(\mathbf{x}^{t}),\label{eq:review-proposed-stepsize-successive}
\end{align}
where $d(\mathbf{x}^{t})$ is the descent defined in (\ref{eq:review-descent-direction}).

The above steps are summarized in Alg. \ref{alg:Successive-approximation-method}.
As a descent algorithm, it generates a monotonically decreasing sequence
$\{h(\mathbf{x}^{t})\}$, and every limit point of $\{\mathbf{x}^{t}\}$
is a stationary point of (\ref{eq:problem-formulation}) (see \cite[Thm. 2]{Yang_ConvexApprox}
for the proof).

\section{\label{sec:Block-SCA}The Proposed Block Successive Convex Approximation
Algorithms}

From a theoretical perspective, Alg. \ref{alg:Successive-approximation-method}
is fully parallelizable. In practice, however, it may not be fully
parallelized when the problem dimension exceeds the hardware's memory
and/or processing capability. We could naively solve the independent
subproblems in Step S1 of Alg. \ref{alg:Successive-approximation-method}
sequentially, for example, in a cyclic order. Once all independent
subproblems are solved, a joint line search is performed as in Step
S2 of Alg. \ref{alg:Successive-approximation-method}.\textbf{ }However,
when the approximation subproblem w.r.t. $\mathbf{x}_{k}$ is being
solved, the solutions of previous approximation subproblems w.r.t.
$(\mathbf{x}_{j})_{j=1}^{k-1}$ are already available, but they are
not exploited.

An alternative is to apply the BCD algorithm, where the variable $\mathbf{x}$
is first divided into blocks $\mathbf{x}=(\mathbf{x}_{k})_{k=1}^{K}$
and the block variables are updated sequentially. Suppose $\mathbf{x}_{k}$
is being updated at iteration $t$, the following optimization problem
w.r.t. the block variable $\mathbf{x}_{k}$ (rather than the full
variable $\mathbf{x}$) is solved while the other block variables
$\mathbf{x}_{-k}\triangleq(\mathbf{x}_{j})_{j\neq k}$ are fixed:
\begin{align}
\mathbf{x}_{k}^{t+1} & =\underset{\mathbf{x}_{k}\in\mathcal{X}_{k}}{\arg\min}\;h(\mathbf{x}_{k},\mathbf{x}_{-k}^{t})\nonumber \\
 & =\underset{\mathbf{x}_{k}\in\mathcal{X}_{k}}{\arg\min}\left\{ f(\mathbf{x}_{k},\mathbf{x}_{-k}^{t})+g_{k}(\mathbf{x}_{k})\right\} .\label{eq:optimization-problem-block-variable}
\end{align}
Convergence to a stationary point of problem (\ref{eq:problem-formulation})
is guaranteed if, for example, $\mathbf{x}_{k}^{t+1}$ is unique \cite{Tseng2001}.
However, the optimization problem in (\ref{eq:optimization-problem-block-variable})
may still not be easy to solve. One approach is to apply Alg. \ref{alg:Successive-approximation-method}
to solve (\ref{eq:optimization-problem-block-variable}) iteratively,
but the resulting algorithm will be of two layers: Alg. \ref{alg:Successive-approximation-method}
keeps iterating in the inner layer until a given accuracy is reached
and the block variable to be updated next is selected in the outer
layer.

To reduce the stringent requirement on the processing capability of
the hardware imposed by the parallel SCA algorithms and the complexity
of the BCD algorithm, we design in this section a BSCA algorithm:
when the block variable $\mathbf{x}_{k}$ is selected at iteration
$t$, all elements of $\mathbf{x}_{k}$ are updated in parallel by
solving an approximation subproblem w.r.t. $\mathbf{x}_{k}$ (rather
than the whole variable $\mathbf{x}$ as in Alg. \ref{alg:Successive-approximation-method})
that is presumably much easier to optimize than the original problem
(\ref{eq:optimization-problem-block-variable}):
\begin{equation}
\mathbb{B}_{k}\mathbf{x}^{t}\triangleq\underset{\mathbf{x}_{k}\in\mathcal{X}_{k}}{\arg\min}\bigl\{\underbrace{\widetilde{f}(\mathbf{x}_{k};\mathbf{x}^{t})+g_{k}(\mathbf{x}_{k})}_{\widetilde{h}(\mathbf{x}_{k};\mathbf{x}^{t})}\bigr\}.\label{eq:hybrid-approximate-problem}
\end{equation}
Note that $\widetilde{f}(\mathbf{x}_{k};\mathbf{x}^{t})$ and $\widetilde{h}(\mathbf{x}_{k};\mathbf{x}^{t})$
defined in (\ref{eq:hybrid-approximate-problem}) is an approximation
function of $f(\mathbf{x}_{k},\mathbf{x}_{-k}^{t})$ and $h(\mathbf{x}_{k},\mathbf{x}_{-k}^{t})$
at a given point $\mathbf{x=x}^{t}$, respectively. We assume that
the approximation function $\widetilde{f}(\mathbf{x};\mathbf{y})$
satisfies the following technical conditions:

\noindent (A1) The function $\widetilde{f}(\mathbf{x}_{k};\mathbf{y})$
is strictly convex in $\mathbf{x}_{k}$ for any given $\mathbf{y}\in\mathcal{X}$;

\noindent (A2) The function $\widetilde{f}(\mathbf{x}_{k};\mathbf{y})$
is continuously differentiable in $\mathbf{x}_{k}$ for any given
$\mathbf{y}\in\mathcal{X}$ and continuous in $\mathbf{y}$ for any
$\mathbf{x}_{k}\in\mathcal{X}$;

\noindent (A3) The gradient of $\widetilde{f}(\mathbf{x}_{k};\mathbf{y})$
and the gradient of $f(\mathbf{x})$ w.r.t. $\mathbf{x}_{k}$ are
identical at $\mathbf{x}=\mathbf{y}$ for any $\mathbf{y}\in\mathcal{X}$,
i.e., $\nabla_{\mathbf{x}_{k}}\widetilde{f}(\mathbf{y}_{k};\mathbf{y})=\nabla_{\mathbf{x}_{k}}f(\mathbf{y})$;

\noindent (A4) A solution $\mathbb{B}_{k}\mathbf{x}^{t}$ exists
for any $\mathbf{x}^{t}\in\mathcal{X}$.

Since the objective function in (\ref{eq:hybrid-approximate-problem})
is strictly convex, $\mathbb{B}_{k}\mathbf{x}^{t}$ is unique. If
$\mathbb{B}_{k}\mathbf{x}^{t}=\mathbf{x}_{k}^{t}$, then $\mathbf{x}_{k}^{t}$
is the optimal point of the optimization problem in (\ref{eq:optimization-problem-block-variable})
given fixed $(\mathbf{x}_{j})_{j\neq k}$ \cite[Prop. 1]{Yang_ConvexApprox}.
We thus consider the case that $\mathbb{B}_{k}\mathbf{x}^{t}\neq\mathbf{x}_{k}^{t}$,
and this implies that
\begin{align}
\widetilde{f}(\mathbb{B}_{k}\mathbf{x}^{t};\mathbf{x}^{t})+g_{k}(\mathbb{B}_{k}\mathbf{x}^{t}) & =\widetilde{h}(\mathbb{B}_{k}\mathbf{x}^{t};\mathbf{x}^{t})\nonumber \\
 & <\widetilde{h}(\mathbf{x}_{k}^{t};\mathbf{x}^{t})=\widetilde{f}(\mathbf{x}_{k}^{t};\mathbf{x}^{t})+g_{k}(\mathbf{x}_{k}^{t}).\label{eq:descent-direction-1}
\end{align}
It follows from the strict convexity of $\widetilde{f}$ and Assumption
(A3) that
\begin{align}
\widetilde{f}(\mathbb{B}_{k}\mathbf{x}^{t};\mathbf{x}^{t})-\widetilde{f}(\mathbf{x}_{k}^{t};\mathbf{x}^{t}) & >(\mathbb{B}_{k}\mathbf{x}^{t}-\mathbf{x}_{k}^{t})^{T}\nabla_{\mathbf{x}_{k}}\widetilde{f}(\mathbf{x}_{k}^{t};\mathbf{x}^{t})\nonumber \\
 & =(\mathbb{B}_{k}\mathbf{x}^{t}-\mathbf{x}_{k}^{t})^{T}\nabla_{\mathbf{x}_{k}}f(\mathbf{x}^{t}).\label{eq:descent-direction-2}
\end{align}
Combining (\ref{eq:descent-direction-1}) and (\ref{eq:descent-direction-2}),
we readily see that $\mathbb{B}_{k}\mathbf{x}^{t}-\mathbf{x}_{k}^{t}$
is a descent direction of $h(\mathbf{x})$ at $\mathbf{x}=\mathbf{x}^{t}$
along the coordinate of $\mathbf{x}_{k}$ in the sense that:
\begin{equation}
d_{k}(\mathbf{x}^{t})\triangleq(\mathbb{B}_{k}\mathbf{x}^{t}-\mathbf{x}_{k}^{t})^{T}\nabla_{\mathbf{x}_{k}}f(\mathbf{x}^{t})+g_{k}(\mathbb{B}_{k}\mathbf{x}^{t})-g_{k}(\mathbf{x}_{k}^{t})<0.\label{eq:descent-direction}
\end{equation}
Then $\mathbf{x}$ is updated according to the following expression:
$\mathbf{x}^{t+1}=(\mathbf{x}_{j}^{t+1})_{j=1}^{K}$ and
\begin{align}
\mathbf{x}_{j}^{t+1} & =\begin{cases}
\mathbf{x}_{k}^{t}+\gamma^{t}(\mathbb{B}_{k}\mathbf{x}^{t}-\mathbf{x}_{k}^{t}), & \textrm{if }j=k,\\
\mathbf{x}_{j}^{t}, & \textrm{otherwise}.
\end{cases}\label{eq:hybrid-variable-update-x}
\end{align}
In other words, only the block variable $\mathbf{x}_{k}$ is updated
while other block variables $(\mathbf{x}_{j})_{j\neq k}$ are equal
to their value at the previous iteration. The stepsize $\gamma^{t}$
in (\ref{eq:hybrid-variable-update-x}) can be determined along the
coordinate of $\mathbf{x}_{k}$ efficiently by the line search introduced
in the previous section, namely, either the exact line search
\begin{equation}
\gamma^{t}=\underset{0\leq\gamma\leq1}{\arg\min}\left\{ \begin{array}{l}
f(\mathbf{x}_{k}^{t}+\gamma(\mathbb{B}_{k}\mathbf{x}^{t}-\mathbf{x}_{k}^{t}),\mathbf{x}_{-k}^{t})\smallskip\\
+\gamma(g_{k}(\mathbb{B}_{k}\mathbf{x}^{t})-g_{k}(\mathbf{x}_{k}^{t}))
\end{array}\right\} ,\label{eq:hybrid-exact-line-search-proposed}
\end{equation}
or the successive line search if the nonconvex differentiable function
in (\ref{eq:hybrid-exact-line-search-proposed}) is still difficult
to optimize: given predefined constants $\alpha\in(0,1)$ and $\beta\in(0,1)$,
the stepsize is set to $\gamma^{t}=\beta^{m_{t}}$, where $m_{t}$
is the smallest nonnegative integer satisfying the inequality:

\begin{align}
 & f(\mathbf{x}_{k}^{t}+\beta^{m}(\mathbb{B}_{k}\mathbf{x}^{t}-\mathbf{x}_{k}^{t}),\mathbf{x}_{-k}^{t})+g_{k}(\mathbf{x}_{k}^{t})+\nonumber \\
 & \quad\beta^{m}(g_{k}(\mathbb{B}_{k}\mathbf{x}^{t})-g_{k}(\mathbf{x}_{k}^{t}))\leq f(\mathbf{x}^{t})+g_{k}(\mathbf{x}_{k}^{t})+\alpha\beta^{m}d_{k}(\mathbf{x}^{t}),\label{eq:hybrid-successive-line-search-proposed}
\end{align}
where $d_{k}(\mathbf{x}^{t})$ is the descent in (\ref{eq:descent-direction}).
Note that the line search in (\ref{eq:hybrid-exact-line-search-proposed})-(\ref{eq:hybrid-successive-line-search-proposed})
is performed along the coordinate of $\mathbf{x}_{k}$ only.

At the next iteration $t+1$, a new block variable is selected and
updated. We consider two commonly used rules to select the block variable,
namely, the cyclic update rule and the random update rule. Note that
both of them are well-known (see \cite{Razaviyayn2014,Scutari_hybrid}),
but we give their definitions for the sake of reference in later developments.

\emph{Cyclic update rule:} The block variables are updated in a cyclic
order. That is, we select the block variable with index\begin{subequations}\label{eq:update-rule}
\begin{equation}
k=\textrm{mod}(t,K)+1.\label{eq:update-rule-cyclic}
\end{equation}

\emph{Random update rule:} The block variables are selected randomly
according to
\begin{equation}
\textrm{Prob}(\mathbf{x}_{k}\textrm{ is updated at iteration }t)=p_{k}^{t}\geq p_{\min}>0,\forall k,\label{eq:update-rule-random}
\end{equation}
\end{subequations}and $\sum_{k}p_{k}^{t}=1$. Any block variable
can be selected with a nonzero probability, and some examples are
given in \cite{Scutari_hybrid}.

The proposed BSCA algorithm is summarized in Alg. \ref{alg:block-Successive-approximation-method},
and its convergence properties are given in the following theorem.

\begin{algorithm}[t]
\textbf{Initialization: }$t=0$, $\mathbf{x}^{0}\in\mathcal{X}$ (arbitrary
but fixed).

Repeat the following steps until convergence:

\begin{enumerate}

\item[\textbf{S1: }]Select the block variable $\mathbf{x}_{k}$ to
be updated according to (\ref{eq:update-rule}).

\item[\textbf{S2: }]Compute $\mathbb{B}_{k}\mathbf{x}^{t}$ according
to (\ref{eq:hybrid-approximate-problem}).

\item[\textbf{S3: }]Determine the stepsize $\gamma^{t}$ by the exact
line search (\ref{eq:hybrid-exact-line-search-proposed}) or the successive
line search (\ref{eq:hybrid-successive-line-search-proposed}).

\item[\textbf{S4: }]Update $\mathbf{x}^{t+1}$ according to (\ref{eq:hybrid-variable-update-x}).

\item[\textbf{S5: }]$t\leftarrow t+1$ and go to \textbf{S1}.

\end{enumerate}

\caption{\label{alg:block-Successive-approximation-method}The proposed block
successive convex approximation algorithm}
\end{algorithm}

\begin{thm}
\label{thm:block-SCA-convergence}Every limit point of the sequence
$\{\mathbf{x}^{t}\}_{t}$ generated by the BSCA algorithm in Alg.
\ref{alg:block-Successive-approximation-method} is a stationary point
of (\ref{eq:problem-formulation}) (with probability 1 for the random
update).
\end{thm}
\begin{IEEEproof}
See Appendix \ref{sec:Proof-of-Theorem-BSCA}.
\end{IEEEproof}
The existence of a limit point is guaranteed if the constraint set
$\mathcal{X}$ in (\ref{eq:problem-formulation}) is bounded or the
objective function $h$ has a bounded lower level set. A sufficient
condition for the latter is that $h$ is coercive, i.e., $h(\mathbf{x})\rightarrow\infty$
as $\left\Vert \mathbf{x}\right\Vert \rightarrow\infty$.

If $\widetilde{f}(\mathbf{x}_{k};\mathbf{x}^{t})$ is a global upper
bound of $f(\mathbf{x}_{k},\mathbf{x}_{-k}^{t})$, we can simply use
a constant unit stepsize $\gamma^{t}=1$, that is,
\begin{equation}
\mathbf{x}_{k}^{t+1}=\mathbb{B}_{k}\mathbf{x}^{t}=\underset{\mathbf{x}_{k}\in\mathcal{X}_{k}}{\arg\min}\;\widetilde{f}(\mathbf{x}_{k};\mathbf{x}^{t})+g_{k}(\mathbf{x}_{k}),\label{eq:variable-update-upper-bound}
\end{equation}
for the reason that the constant unit stepsize always yields a larger
decrease than the successive line search and the convergence is guaranteed
(see the discussion on Assumption (A6) in \cite[Sec. III]{Yang_ConvexApprox}).
In this case, update (\ref{eq:variable-update-upper-bound}) has the
same form as BSUM \cite{Razaviyayn2013}. However, their convergence
conditions and techniques are different and do not imply each other.
As a matter of fact, stronger results may be obtained, see \cite{Xu2017}.

There are several commonly used choices of approximation function
$\widetilde{f}(\mathbf{x}_{k};\mathbf{x}^{t})$, for example, the
linear approximation and the quadratic approximation. We refer to
\cite[Sec. III-B]{Yang_ConvexApprox} and \cite[Sec. II.2.1]{Scutari2018}
for more details and just comment on the following important cases.

\emph{Quadratic approximation:}
\begin{equation}
\widetilde{f}(\mathbf{x}_{k};\mathbf{x}^{t})=(\mathbf{x}_{k}-\mathbf{x}_{k}^{t})^{T}\nabla_{k}f(\mathbf{x}^{t})+\frac{c_{k}^{t}}{2}\left\Vert \mathbf{x}_{k}-\mathbf{x}_{k}^{t}\right\Vert _{2}^{2},\label{eq:approximation-quadratic}
\end{equation}
where $c_{k}^{t}$ is a positive scalar. If $\nabla_{k}f(\mathbf{x}_{k},\mathbf{x}_{-k}^{t})$
is Lipschitz continuous, $\widetilde{f}(\mathbf{x}_{k};\mathbf{x}^{t})$
would be a global upper bound of $f(\mathbf{x}_{k},\mathbf{x}_{-k}^{t})$
when $c_{k}^{t}$ is sufficiently large. In this case, the variable
update reduces to the well-known proximal operator:
\[
\mathbf{x}_{k}^{t+1}=\underset{\mathbf{x}_{k}}{\arg\min}\left\{ \begin{array}{l}
(\mathbf{x}_{k}-\mathbf{x}_{k}^{t})^{T}\nabla_{k}f(\mathbf{x}^{t})\smallskip\\
+\frac{c_{k}^{t}}{2}\left\Vert \mathbf{x}_{k}-\mathbf{x}_{k}^{t}\right\Vert _{2}^{2}+g_{k}(\mathbf{x}_{k})
\end{array}\right\} ,
\]
and this is also known as the proximal linear approximation. If we
incorporate a stepsize as in the proposed BSCA algorithm, $\widetilde{f}(\mathbf{x}_{k};\mathbf{x}^{t})$
is strictly convex as long as $c_{k}^{t}$ is positive and the convergence
is thus guaranteed by Theorem \ref{thm:block-SCA-convergence} (even
when $\nabla_{k}f(\mathbf{x}_{k},\mathbf{x}_{-k}^{t})$ is not Lipschitz
continuous).

\emph{Best-response approximation \#1: }If $f(\mathbf{x})$ is strictly
convex in each element of the block variable $\mathbf{x}_{k}=(x_{i_{k}})_{i_{k}=1}^{I_{k}}$,
the ``best-response'' type approximation function is
\begin{equation}
\widetilde{f}(\mathbf{x}_{k};\mathbf{x}^{t})=\sum_{i_{k=1}}^{I_{k}}f(x_{i_{k}},(x_{j_{k}}^{t})_{j_{k}\neq i_{k}},\mathbf{x}_{-k}^{t}),\label{eq:approximation-best-response-1}
\end{equation}
and it is \emph{not} a global upper bound of $f(\mathbf{x}_{k},\mathbf{x}_{-k}^{t})$.
Note that $f(\mathbf{x}_{k},\mathbf{x}_{-k}^{t})$ is not necessarily
convex in $\mathbf{x}_{k}$ and the best-response approximation is
different from the above proximal linear approximation and thus cannot
be obtained from existing algorithmic frameworks \cite{Xu2013a,Xu2017,Bolte2014}.

\emph{Best-response approximation \#2:} If $f(\mathbf{x})$ is furthermore
strictly convex in $\mathbf{x}_{k}$, an alternative ``best-response''
type approximation function is
\begin{equation}
\widetilde{f}(\mathbf{x}_{k};\mathbf{x}^{t})=f(\mathbf{x}_{k},\mathbf{x}_{-k}^{t}).\label{eq:approximation-best-response-2}
\end{equation}
The approximation function in (\ref{eq:approximation-best-response-2})
is a trivial upper bound of $f(\mathbf{x}_{k},\mathbf{x}_{-k}^{t})$,
and the BSCA algorithm (\ref{eq:variable-update-upper-bound}) reduces
to the BCD algorithm (\ref{eq:optimization-problem-block-variable}).
Adopting the approximation in (\ref{eq:approximation-best-response-2})
usually leads to fewer iterations than (\ref{eq:approximation-best-response-1}),
as (\ref{eq:approximation-best-response-2}) is a ``better'' approximation
in the sense that it is on the basis of the block variable $\mathbf{x}_{k}$,
while the approximation in (\ref{eq:approximation-best-response-1})
is on the basis of each element of $\mathbf{x}_{k}$, namely, $x_{i_{k}}$
for all $i_{k}=1,\ldots,I_{k}$. Nevertheless, the approximation function
(\ref{eq:approximation-best-response-1}) may be easier to optimize
than (\ref{eq:approximation-best-response-2}) as the component functions
are separable and each component function is a scalar function. This
reflects the universal tradeoff between the number of iterations and
the complexity per iteration.

\emph{Partial linearization approximation:} Consider the function
$f=f_{1}(f_{2}(\mathbf{x}))$ where $f_{1}$ is smooth and convex
and $f_{2}$ is smooth. We can adopt the ``partial linearization''
approximation where $f_{2}(\mathbf{x})$ is linearized while $f_{1}$
is left unchanged:
\begin{align}
\widetilde{f}(\mathbf{x}_{k};\mathbf{x}^{t})=\; & f_{1}(f_{2}(\mathbf{x}^{t})+(\mathbf{x}_{k}-\mathbf{x}_{k}^{t})\nabla_{k}f_{2}(\mathbf{x}))\nonumber \\
 & +\frac{c_{k}^{t}}{2}\left\Vert \mathbf{x}_{k}-\mathbf{x}_{k}^{t}\right\Vert _{2}^{2},\label{eq:approximation-partial-linearization}
\end{align}
where $c_{k}^{t}$ is a positive scalar. The quadratic regularization
is incorporated to make the approximation function strictly convex.
It can be verified by using the chain rule that
\[
\nabla_{k}\widetilde{f}(\mathbf{x}_{k}^{t};\mathbf{x}^{t})=\nabla f_{1}(f_{2}(\mathbf{x}^{t}))\nabla_{k}f_{2}(\mathbf{x}^{t})=\nabla f(\mathbf{x}^{t}).
\]
The partial linearization approximation is expected to yield faster
convergence than quadratic approximation because the convexity of
function $f_{1}$ is preserved in (\ref{eq:approximation-partial-linearization}).

\emph{Hybrid approximation:} For the above composition function $f=f_{1}(f_{2}(\mathbf{x}))$,
we can also adopt a hybrid approximation by further approximating
the partial linearization approximation (\ref{eq:approximation-partial-linearization})
by the best-response approximation (\ref{eq:approximation-best-response-1}):
\begin{align}
\widetilde{f}(\mathbf{x}_{k};\mathbf{x}^{t})=\; & \sum_{i_{k}=1}^{I_{k}}f_{1}(f_{2}(\mathbf{x}^{t})+(x_{i_{k}}-x_{i_{k}}^{t})\nabla_{i_{k}}f_{2}(\mathbf{x}))\nonumber \\
 & +\frac{c_{k}^{t}}{2}\left\Vert \mathbf{x}_{k}-\mathbf{x}_{k}^{t}\right\Vert _{2}^{2}.\label{eq:approximation-hybrid-linearization}
\end{align}
The hybrid approximation function (\ref{eq:approximation-hybrid-linearization})
is separable among the elements of $\mathbf{x}_{k}$, while it is
not necessarily the case for the partial linearization approximation
(\ref{eq:approximation-partial-linearization}). The separable structure
is desirable when $g_{k}(\mathbf{x}_{k})$ is also separable among
the elements of $\mathbf{x}_{k}$ (for example $g_{k}(\mathbf{x}_{k})=\left\Vert \mathbf{x}_{k}\right\Vert _{1}$),
because the approximation subproblem (\ref{eq:hybrid-approximate-problem})
would further boil down to parallel scalar problems. We remark that
the partial linearization approximation and the hybrid approximation
are only foreseen by SCA framework and cannot be obtained from other
existing algorithmic frameworks \cite{Xu2013a,Xu2017,Bolte2014}.
\begin{rem}
\label{rem:flexibility}The above approximation is on the basis of
blocks and it may be different from block to block. For example, consider
$f(\mathbf{x})=f_{1}(\mathbf{x}_{1},\mathbf{x}_{2},f_{2}(\mathbf{x}_{3}))$
where $f_{1}$ is strictly convex in $\mathbf{x}_{1}$, nonconvex
in $\mathbf{x}_{2}$, and convex in $f_{2}(\mathbf{x}_{3})$. Then
we can adopt the best-response approximation for $\mathbf{x}_{1}$,
the quadratic approximation for $\mathbf{x}_{2}$, and partial linearization
(or hybrid) approximation for $\mathbf{x}_{3}$:
\begin{align*}
\widetilde{f}(\mathbf{x}_{1};\mathbf{x}^{t})=\; & f_{1}(\mathbf{x}_{1},\mathbf{x}_{2}^{t},\mathbf{x}_{3}^{t}),\\
\widetilde{f}(\mathbf{x}_{2};\mathbf{x}^{t})=\; & (\mathbf{x}_{2}-\mathbf{x}_{2}^{t})^{T}\nabla_{\mathbf{x}_{2}}f_{1}(\mathbf{x}^{t})+\frac{\tau^{t}}{2}\left\Vert \mathbf{x}_{2}-\mathbf{x}_{2}^{t}\right\Vert _{2}^{2},\\
\widetilde{f}(\mathbf{x}_{3};\mathbf{x}^{t})=\; & f_{1}(\mathbf{x}_{1}^{t},\mathbf{x}_{2}^{t},f_{2}(\mathbf{x}_{3}^{t})+(\mathbf{x}_{3}-\mathbf{x}_{3}^{t})^{T}\nabla_{\mathbf{x}_{3}}f_{2}(\mathbf{x}_{3}^{t}))\\
 & +\frac{c^{t}}{2}\left\Vert \mathbf{x}_{3}-\mathbf{x}_{3}^{t}\right\Vert _{2}^{2}.
\end{align*}
The most suitable approximation always depends on the application
and the universal tradeoff between the number of iterations and the
complexity per iteration. SCA offers sufficient flexibility to address
this tradeoff.
\end{rem}
The proposed BSCA algorithm described in Alg. \ref{alg:block-Successive-approximation-method}
is complementary to the parallel SCA algorithm in Alg. \ref{alg:Successive-approximation-method}.
On the one hand, the update of the elements of a particular block
variable in the BSCA algorithm is based on the same principle as in
the parallel SCA algorithm, namely, to obtain the descent direction
by minimizing a convex approximation function and to calculate the
stepsize by the line search scheme. On the other hand, in contrast
to the parallel update in the parallel SCA algorithm, the block variables
are updated sequentially in the BSCA algorithm, and it poses a less
demanding requirement on the memory/processing unit.

We draw some comments on the proposed BSCA algorithm.

\textbf{On the connection to traditional BCD algorithms.} The point
$\mathbf{x}_{k}^{t+1}$ in (\ref{eq:hybrid-variable-update-x}) is
obtained by moving from the current point $\mathbf{x}_{k}^{t}$ along
a descent direction $\mathbb{B}_{k}\mathbf{x}^{t}-\mathbf{x}_{k}^{t}$.
On the one hand, $\mathbf{x}_{k}^{t+1}$ is in general not the best-response
employed in the traditional BCD algorithm (\ref{eq:optimization-problem-block-variable}).
That is,
\[
f(\mathbf{x}_{k}^{t+1},\mathbf{x}_{-k}^{t})+g_{k}(\mathbf{x}_{k})\geq\min_{\mathbf{x}_{k}\in\mathcal{X}_{k}}f(\mathbf{x}_{k},\mathbf{x}_{-k}^{t})+g_{k}(\mathbf{x}_{k}).
\]
Therefore, the proposed algorithm is essentially an inexact BCD algorithm.
On the other hand, Theorem \ref{thm:block-SCA-convergence} establishes
that eventually there is no loss of optimality adopting inexact solutions
as long as the approximate functions satisfy the assumptions (A1)-(A4).

\textbf{On the flexibility.} The assumptions (A1)-(A4) on the approximation
function are quite general and they include many existing algorithms
as a special case (see Remark \ref{rem:flexibility}). The proposed
approximation function does not have to be a global upper bound of
the original function, but a stepsize is needed to avoid aggressive
update.

\textbf{On the convergence speed.} The mild assumptions on the approximation
functions allow us to design an approximation function that exploits
the original problem's structure (such as the partial convexity in
(\ref{eq:approximation-best-response-1})-(\ref{eq:approximation-best-response-2}))
and this leads to faster convergence. The use of line search also
attributes to a faster convergence than decreasing stepsizes used
in literature, for example \cite{Scutari_BigData,Razaviyayn2014}.

\textbf{On the complexity.} The proposed BSCA exhibits low complexity
for several reasons. Firstly, the exact line search consists of minimizing
a differentiable function. Although this incurs additional complexity
compared with pre-determined stepsizes, in many signal processing
and machine learning applications, the line search admits a closed-form
solution, as we shown later in the example applications. In the successive
line search, the nonsmooth function $g$ only needs to be evaluated
once at the point $\mathbb{B}_{k}\mathbf{x}^{t}$. Secondly, the problem
size that can be handled by the BSCA algorithm is much larger.

\textbf{On the convergence conditions.} The strict convexity of the
approximation function $\tilde{f}(\mathbf{x}_{k},\mathbf{y})$ according
to (A1) is stronger than the convexity assumption in the approximation
function for parallel SCA (reviewed in Sec. \ref{sec:Review-of-SCA}).
This is to guarantee the approximation subproblem has a unique solution,
which is essential to ensure the convergence of the block update.
The subsequence convergence of the BSCA algorithm is established under
fairly weak assumptions in Theorem \ref{thm:block-SCA-convergence}.
Compared with \cite{Razaviyayn2013,Scutari_BigData}, the BSCA algorithm
is applicable for nonsmooth nonconvex optimization problems, and it
converges even when the gradient of $f$ is not Lipschitz continuous,
respectively. Nevertheless, the subsequence convergence is weaker
than the sequence convergence \cite{Xu2013a,Bolte2014,Chouzenoux2016,Xu2017,Bolte2018}:
in theory it is possible that two convergent subsequences converge
to different stationary points, so the whole sequence may diverge.

\section{\label{sec:Inexact-Block-SCA}The Proposed Inexact Block Successive
Convex Approximation Algorithms}

In the previous section, the approximation subproblem in (\ref{eq:hybrid-approximate-problem})
is assumed to be solved exactly, and this assumption is satisfied
when $\mathbb{B}_{k}\mathbf{x}^{t}$ has a closed-form expression.
However, if $\mathbb{B}_{k}\mathbf{x}^{t}$ does not have a closed-form
expression for some choice of the approximation function, it must
be found numerically by an iterative algorithm. In this case, Alg.
\ref{alg:block-Successive-approximation-method} would consist of
two layers: the outer layer with index $t$ follows the same procedure
as Alg. \ref{alg:block-Successive-approximation-method}, while the
inner layer comprising the iterative algorithm for (\ref{eq:hybrid-approximate-problem})
is nested under S2 of Alg. \ref{alg:block-Successive-approximation-method}.
As most iterative algorithms exhibit asymptotic convergence only,
in practice, they are terminated when we obtain an approximate solution,
denoted as $\widetilde{\mathbf{x}}_{k}^{t}$, which is ``sufficiently
accurate'' in the sense that the so-called error bound $\left\Vert \widetilde{\mathbf{x}}_{k}^{t}-\mathbb{B}_{k}\mathbf{x}^{t}\right\Vert \leq\epsilon^{t}$
for some small $\epsilon^{t}$ that decreases to zero as $t$ increases
\cite{Scutari_BigData,Scutari_hybrid}. Nevertheless, results on the
error bound are mostly available when $\widetilde{f}(\mathbf{x}_{k};\mathbf{x}^{t})$
is strongly convex, see \cite[Ch. 6]{Facchinei&Pang}. In general,
the error bound is very difficult to verify in practice.

In this section, we propose an inexact BSCA algorithm where (\ref{eq:hybrid-approximate-problem})
is solved inexactly, but we do not pose any quantitative requirement
on its error bound. The central idea is in (\ref{eq:review-descent-direction}):
replacing $\mathbb{B}_{k}\mathbf{x}^{t}$ by any point $\widetilde{\mathbf{x}}_{k}^{t}$
and repeating the same steps in (\ref{eq:descent-direction-1})-(\ref{eq:descent-direction}),
we see that $\widetilde{\mathbf{x}}_{k}^{t}-\mathbf{x}_{k}^{t}$ is
a descent direction of $h(\mathbf{x}_{k},\mathbf{x}_{-k}^{t})$ at
$\mathbf{x}_{k}=\mathbf{x}_{k}^{t}$ if
\[
\widetilde{f}(\widetilde{\mathbf{x}}_{k}^{t};\mathbf{x}^{t})+g_{k}(\widetilde{\mathbf{x}}_{k}^{t})-(\widetilde{f}(\mathbf{x}_{k}^{t};\mathbf{x}^{t})+g_{k}(\mathbf{x}_{k}^{t}))<0.
\]
Such a point $\widetilde{\mathbf{x}}_{k}^{t}$ can be obtained by
running the standard parallel SCA algorithm (reviewed in Sec. \ref{sec:Review-of-SCA})
to solve the approximation subproblem (\ref{eq:hybrid-approximate-problem})
in S2 of Alg. \ref{alg:block-Successive-approximation-method} for
a\emph{ finite} number of iterations only. At iteration $\tau$ of
the inner layer, we define $\widetilde{f}^{i}(\mathbf{x}_{k};\mathbf{x}_{k}^{t,\tau},\mathbf{x}^{t})$
(the superscript ``\emph{i}'' stands for inner) as an approximation
of the (outer-layer) approximation function $\widetilde{f}(\mathbf{x}_{k};\mathbf{x}^{t})$
at $\mathbf{x}_{k}=\mathbf{x}_{k}^{t,\tau}$ (with $\mathbf{x}_{k}^{t,0}=\mathbf{x}_{k}^{t}$),
and solve the inner-layer approximation subproblem
\begin{equation}
\mathbb{B}_{k}\overline{\mathbf{x}}^{t,\tau}\in\underset{\mathbf{x}_{k}\in\mathcal{X}_{k}}{\arg\min}\;\{\underbrace{\widetilde{f}^{i}(\mathbf{x}_{k};\mathbf{x}_{k}^{t,\tau},\mathbf{x}^{t})+g_{k}(\mathbf{x}_{k})}_{\widetilde{h}^{i}(\mathbf{x}_{k};\mathbf{x}_{k}^{t,\tau},\mathbf{x}^{t})}\},\label{eq:inexact-approximate-problem}
\end{equation}
where $\overline{\mathbf{x}}^{t,\tau}\triangleq(\mathbf{x}_{k}^{t,\tau},\mathbf{x}^{t})$.
Presumably this problem is designed to be much easier to solve exactly
than the outer-layer approximation subproblem in (\ref{eq:hybrid-approximate-problem}),
for example, a closed-form solution exists; such an example will be
given later in Sec. \ref{sec:Example-applications}-C. We assume the
inner-layer approximation function $\widetilde{f}^{i}$ satisfies
the following technical assumptions.

\noindent (B1) The function $\widetilde{f}_{k}^{i}(\mathbf{x}_{k};\mathbf{x}_{k}^{t,\tau},\mathbf{x}^{t})$
is convex in $\mathbf{x}_{k}$ for any given $\mathbf{x}_{k}^{t,\tau}\in\mathcal{X}_{k}$
and $\mathbf{x}^{t}\in\mathcal{X}$;

\noindent (B2) The function $\widetilde{f}_{k}^{i}(\mathbf{x}_{k};\mathbf{x}_{k}^{t,\tau},\mathbf{x}^{t})$
is continuously differentiable in $\mathbf{x}_{k}$ for any given
$\mathbf{x}_{k}^{t,\tau}\in\mathcal{X}_{k}$ and $\mathbf{x}^{t}\in\mathcal{X}$,
and continuous in $\mathbf{x}_{k}^{t,\tau}$ and $\mathbf{x}^{t}$
for any $\mathbf{x}_{k}\in\mathcal{X}_{k}$;

\noindent (B3) The gradient of $\widetilde{f}_{k}^{i}(\mathbf{x}_{k};\mathbf{x}_{k}^{t,\tau},\mathbf{x}^{t})$
and the gradient of $\widetilde{f}(\mathbf{x}_{k};\mathbf{x}^{t})$
w.r.t. $\mathbf{x}_{k}$ are identical at $\mathbf{x}_{k}=\mathbf{x}_{k}^{t,\tau}$
for any $\mathbf{x}^{t}\in\mathcal{X}$, i.e., $\nabla_{\mathbf{x}_{k}}\widetilde{f}^{i}(\mathbf{x}_{k}^{t,\tau};\mathbf{x}_{k}^{t,\tau},\mathbf{x}^{t})=\nabla_{\mathbf{x}_{k}}\widetilde{f}(\mathbf{x}_{k}^{t,\tau};\mathbf{x}^{t})$;

\noindent (B4) A solution $\mathbb{B}_{k}\overline{\mathbf{x}}^{t,\tau}$
exists for any $(\mathbf{x}_{k}^{t,\tau},\mathbf{x}^{t})$;

\noindent (B5) For any bounded sequence $\{\mathbf{x}_{k}^{t,\tau}\}_{\tau}$,
the sequence $\{\mathbb{B}_{k}\overline{\mathbf{x}}^{t,\tau}\}$ is
also bounded.

Since $\widetilde{h}^{i}(\mathbf{x}_{k};\mathbf{x}_{k}^{t,\tau},\mathbf{x}^{t})$
is convex in $\mathbf{x}_{k}$ and its (global) minimum value is achieved
at $\mathbb{B}_{k}\overline{\mathbf{x}}^{t,\tau}$, it is either
\begin{equation}
\widetilde{h}^{i}(\mathbb{B}_{k}\overline{\mathbf{x}}^{t,\tau};\mathbf{x}_{k}^{t,\tau},\mathbf{x}_{-k}^{t})=\widetilde{h}^{i}(\mathbf{x}_{k}^{t,\tau};\mathbf{x}_{k}^{t,\tau},\mathbf{x}_{-k}^{t})\label{eq:possibility-1}
\end{equation}
or
\begin{equation}
\widetilde{h}^{i}(\mathbb{B}_{k}\overline{\mathbf{x}}^{t,\tau};\mathbf{x}_{k}^{t,\tau},\mathbf{x}_{-k}^{t})<\widetilde{h}^{i}(\mathbf{x}_{k}^{t,\tau};\mathbf{x}_{k}^{t,\tau},\mathbf{x}_{-k}^{t}).\label{eq:possibility-2}
\end{equation}
On the one hand, if (\ref{eq:possibility-1}) is true, $\mathbf{x}_{k}^{t,\tau}=\mathbb{B}_{k}\mathbf{x}^{t}$
and the outer-layer approximation subproblem in (\ref{eq:hybrid-approximate-problem})
has been solved exactly \cite[Prop. 1]{Yang_ConvexApprox}. On the
other hand, (\ref{eq:possibility-2}) implies that $\mathbb{B}_{k}\overline{\mathbf{x}}^{t,\tau}-\mathbf{x}_{k}^{t,\tau}$
is a descent direction of the outer-layer approximation function $\widetilde{h}(\mathbf{x}_{k};\mathbf{x}^{t})$
at $\mathbf{x}_{k}=\mathbf{x}_{k}^{t,\tau}$, i.e.,
\begin{align}
d_{k}(\overline{\mathbf{x}}^{t,\tau})\triangleq\; & \nabla\widetilde{f}(\mathbf{x}_{k}^{t,\tau};\mathbf{x}^{t})^{T}(\mathbb{B}_{k}\overline{\mathbf{x}}^{t,\tau}-\mathbf{x}_{k}^{t,\tau})\nonumber \\
 & +g_{k}(\mathbb{B}_{k}\overline{\mathbf{x}}^{t,\tau})-g_{k}(\mathbf{x}_{k}^{t,\tau})<0.\label{eq:inexact-descent}
\end{align}
Therefore we can update $\mathbf{x}_{k}^{t,\tau+1}$ by
\begin{equation}
\mathbf{x}_{k}^{t,\tau+1}=\mathbf{x}_{k}^{t,\tau}+\gamma^{t,\tau}(\mathbb{B}_{k}\overline{\mathbf{x}}^{t,\tau}-\mathbf{x}_{k}^{t}),\label{eq:inexact-variable-update}
\end{equation}
where $\gamma^{t,\tau}$ is calculated by either the exact line search
along the coordinate of $\mathbf{x}_{k}$ over the outer-layer approximation
function $\widetilde{h}(\mathbf{x}_{k};\mathbf{x}^{t})$ (rather than
the original function $h(\mathbf{x})$):
\begin{equation}
\gamma^{t,\tau}=\underset{0\leq\gamma\leq1}{\arg\min}\left\{ \begin{array}{l}
\widetilde{f}(\mathbf{x}_{k}^{t,\tau}+\gamma(\mathbb{B}_{k}\overline{\mathbf{x}}^{t,\tau}-\mathbf{x}_{k}^{t,\tau});\mathbf{x}^{t})\smallskip\\
+\gamma(g_{k}(\mathbb{B}_{k}\overline{\mathbf{x}}^{t,\tau})-g_{k}(\mathbf{x}_{k}^{t,\tau}))
\end{array}\right\} ,\label{eq:inexact-line-search-exact}
\end{equation}
or the successive line search: given predefined constants $\alpha\in(0,1)$
and $\beta\in(0,1)$, the stepsize is set to $\gamma^{t,\tau}=\beta^{m_{t,\tau}}$,
where $m_{t,\tau}$ is the smallest nonnegative integer $m$ satisfying
\begin{align}
 & \widetilde{f}(\mathbf{x}_{k}^{t,\tau}+\beta^{m}(\mathbb{B}_{k}\overline{\mathbf{x}}^{t,\tau}-\mathbf{x}_{k}^{t,\tau});\mathbf{x}^{t})\nonumber \\
\leq\; & \widetilde{f}(\mathbf{x}_{k}^{t,\tau};\mathbf{x}^{t})+\beta^{m}(\alpha d_{k}(\overline{\mathbf{x}}^{t,\tau})-(g_{k}(\mathbb{B}_{k}\overline{\mathbf{x}}^{t,\tau})-g_{k}(\mathbf{x}_{k}^{t,\tau}))),\label{eq:inexact-line-search-successive}
\end{align}
where $d_{k}(\mathbf{x}^{t,\tau})$ is the descent defined in (\ref{eq:inexact-descent}).

After repeating the process specified in (\ref{eq:inexact-approximate-problem})-(\ref{eq:inexact-line-search-successive})
for a finite number of iterations denoted by $\bar{\tau}_{t}$, we
set $\widetilde{\mathbf{x}}_{k}^{t}=\mathbf{x}_{k}^{t,\bar{\tau}_{t}}$
and compute the stepsize $\gamma^{t}$ by the line search (\ref{eq:hybrid-exact-line-search-proposed})
or (\ref{eq:hybrid-successive-line-search-proposed}) (therein $\mathbb{B}_{k}\mathbf{x}^{t}$
should be replaced by $\widetilde{\mathbf{x}}_{k}^{t}$). Then $\mathbf{x}^{t+1}=(\mathbf{x}_{j}^{t+1})_{j=1}^{K}$
is set according to
\begin{align}
\mathbf{x}_{j}^{t+1} & =\begin{cases}
\mathbf{x}_{k}^{t}+\gamma^{t}(\widetilde{\mathbf{x}}_{k}^{t}-\mathbf{x}_{k}^{t}), & \textrm{if }j=k,\\
\mathbf{x}_{j}^{t}, & \textrm{otherwise}.
\end{cases}\label{eq:hybrid-inexact-variable-update-x}
\end{align}
The number of inner-layer iterations $\bar{\tau}_{t}$ is a finite
number and may be varying from iteration to iteration. The above procedure
is formally summarized in Alg. \ref{alg:block-Successive-approximation-method-inexact}.

The sequence $\{\widetilde{h}(\mathbf{x}_{k}^{t,\tau};\mathbf{x}^{t})\}_{\tau}$
is monotonically decreasing, but lower bounded by the minimum value
$\widetilde{h}(\mathbb{B}_{k}\mathbf{x}^{t};\mathbf{x}^{t})$:
\begin{align*}
\widetilde{h}(\mathbf{x}_{k}^{t};\mathbf{x}_{-k}^{t})=\widetilde{h}(\mathbf{x}_{k}^{t,0};\mathbf{x}^{t})>\ldots & >\widetilde{h}(\mathbf{x}_{k}^{t,\bar{\tau}_{t}-1};\mathbf{x}^{t})\\
 & >\widetilde{h}(\mathbf{x}_{k}^{t,\bar{\tau}_{t}};\mathbf{x}^{t})\geq\widetilde{h}(\mathbb{B}_{k}\mathbf{x}^{t};\mathbf{x}^{t}).
\end{align*}
This also implies that $\mathbf{x}_{k}^{t,\bar{\tau}_{t}}$ is in
general not an optimal point of the outer-layer approximation subproblem
in (\ref{eq:hybrid-approximate-problem}). However, every limit point
of the sequence $\{\mathbf{x}_{k}^{t,\tau}\}_{\tau}$ is an optimal
point \cite[Thm. 2]{Yang_ConvexApprox}, which is unique in view of
the strict convexity of $\widetilde{h}(\mathbf{x}_{k};\mathbf{x}^{t})$.
Therefore the whole sequence $\{\mathbf{x}_{k}^{t,\tau}\}$ converges
to $\mathbb{B}_{k}\mathbf{x}^{t}$:
\[
\lim_{\tau\rightarrow\infty}\mathbf{x}_{k}^{t,\tau}=\mathbb{B}_{k}\mathbf{x}^{t}.
\]

\begin{algorithm}[t]
\textbf{Initialization: }$t=0$, $\mathbf{x}^{0}\in\mathcal{X}$ (arbitrary
but fixed).

Repeat the following steps until convergence:

\begin{enumerate}

\item[\textbf{S1: }]Select the block variable $\mathbf{x}_{k}$ to
be updated according to the cyclic rule (\ref{eq:update-rule-cyclic})
or the random rule (\ref{eq:update-rule-random}).

\item[\textbf{S2: }]Compute $\widetilde{\mathbf{x}}_{k}^{t}$ by
the following steps:

\textbf{}%
\colorbox{lightgray}{\begin{minipage}[t]{0.88\columnwidth}%
\begin{enumerate}

\item[\textbf{S2.0:}] Set $\tau=0$ and $\mathbf{x}_{k}^{t,0}=\mathbf{x}_{k}^{t}$.

\item[\textbf{S2.1:}] Compute $\mathbb{B}_{k}\mathbf{x}^{t,\tau}$
according to (\ref{eq:inexact-approximate-problem}).

\item[\textbf{S2.3:}] Compute $\gamma^{t,\tau}$ by the line search
(\ref{eq:inexact-line-search-exact}) or (\ref{eq:inexact-line-search-successive}).

\item[\textbf{S2.4:}] Update $\mathbf{x}_{k}^{t,\tau+1}$ according
to (\ref{eq:inexact-variable-update}).

\item[\textbf{S2.2:}] If $\tau+1=\bar{\tau}_{t}$, set $\widetilde{\mathbf{x}}_{k}^{t}=\mathbf{x}_{k}^{t,\bar{\tau}_{t}}$
and go to \textbf{S3}. Otherwise $\tau\leftarrow\tau+1$ and go to
\textbf{S2.1}.

\end{enumerate}%
\end{minipage}}

\item[\textbf{S3: }]Given the update direction $\widetilde{\mathbf{x}}_{k}^{t}-\mathbf{x}_{k}^{t}$,
compute the stepsize $\gamma^{t}$ by the exact line search (\ref{eq:hybrid-exact-line-search-proposed})
or the successive line search (\ref{eq:hybrid-successive-line-search-proposed}).

\item[\textbf{S4: }]Update $\mathbf{x}^{t+1}$ according to (\ref{eq:hybrid-inexact-variable-update-x}).

\item[\textbf{S5: }]$t\leftarrow t+1$ and go to \textbf{S1}.

\end{enumerate}

\caption{\label{alg:block-Successive-approximation-method-inexact}The proposed
inexact block successive convex approximation algorithm}
\end{algorithm}

\begin{thm}
\label{thm:inexact-SCA}If Assumptions (A1)-(A4) and (B1)-(B5) are
satisfied, then every limit point of the sequence $\{\mathbf{x}^{t}\}$
generated by Alg. \ref{alg:block-Successive-approximation-method-inexact}
is a stationary point of Problem (\ref{eq:problem-formulation}) (with
probability 1 for the random update).
\end{thm}
\begin{IEEEproof}
See Appendix \ref{sec:Proof-of-Theorem-inexact-BSCA}.
\end{IEEEproof}
A straightforward choice of $\widetilde{f}^{i}(\mathbf{x}_{k};\mathbf{x}_{k}^{t,\tau},\mathbf{x}_{-k}^{t})$
is
\begin{equation}
\widetilde{f}^{i}(\mathbf{x}_{k};\mathbf{x}_{k}^{t,\tau},\mathbf{x}^{t})=\sum_{k_{i}=1}^{k_{I}}\widetilde{f}(x_{k_{i}},\mathbf{x}_{-k_{i}}^{t,\tau};\mathbf{x}^{t}),\label{eq:best-response-approximation-InnerLayer}
\end{equation}
where $\mathbf{x}_{k}=(x_{k_{i}})_{k_{i}=1}^{k_{I}}$. It is strictly
convex because $\widetilde{f}(\mathbf{x}_{k};\mathbf{x}^{t})$ is
strictly convex in $\mathbf{x}_{k}$ in view of Assumption (A1) and
thus individually strictly convex in each element of $\mathbf{x}_{k}$.

We remark once more that the salient feature of the inexact BSCA algorithm
is that when a (parallel) SCA-based algorithm is applied to solve
the approximation subproblem (\ref{eq:hybrid-approximate-problem}),
it can be terminated after a finite number of iterations without checking
the solution accuracy. Note that the use of a SCA-based algorithm
is not a restrictive assumption as it includes as special cases a
fairly large number of existing algorithms, such as proximal algorithms,
gradient-based algorithms and parallel BCD algorithms \cite[Sec. III-B]{Yang_ConvexApprox}.
Nevertheless, it is not difficult to see that the SCA-based algorithm
nested under Step S2 of Alg. \ref{alg:block-Successive-approximation-method-inexact}
can also be replaced by any other algorithm, as long as it is a closed
mapping\footnote{A mapping $\mathbb{B}:\mathbf{x}\rightarrow\mathbb{B}\mathbf{x}$
is closed if $\mathbf{x}^{t}\rightarrow\mathbf{x}^{\star}$ and $\mathbb{B}\mathbf{x}^{t}\rightarrow\mathbf{y}^{\star}$
for some $\mathbf{x}^{\star}$ and $\mathbf{y}^{\star}$, then $\mathbf{y}^{\star}\in\mathbb{B}\mathbf{x}^{\star}$
\cite{Berge1963}.} that can produce a point $\widetilde{\mathbf{x}}_{k}^{t}$ that has
a lower objective value than $\mathbf{x}_{k}^{t}$. This observation
has profound implications in both theory and practice. From the theoretical
perspective, the complicated error bound in \cite{Scutari_BigData,Scutari_hybrid}
is no longer needed and the convergence condition is significantly
relaxed. Besides, the proposed algorithm extends the inexact SCA algorithm
in \cite{Patriksson1995} where the convergence is proved under the
traditional exact line search in the spirit of (\ref{eq:review-traditional-line-search})
only. From the practical perspective, this leads to extremely easy
implementation without any loss in optimality. We further show through
simulations in Sec. \ref{sec:Example-applications}-B that it is sometimes
not necessary to solve the approximation subproblem (\ref{eq:hybrid-approximate-problem})
with a high precision.

\section{\label{sec:Example-applications}Applications in Sparse Signal Estimation
and Machine Learning}

\subsection{Joint Estimation of Low-Rank and Sparse Signals}

 \newcounter{MYtempeqncnt1} \begin{figure*}[t] \normalsize \setcounter{MYtempeqncnt1}{\value{equation}} \setcounter{equation}{36} \vspace*{4pt}
\begin{equation}
\gamma^{t}=\left[-\frac{\textrm{tr}((\mathbf{P}^{t}\mathbf{Q}^{t}+\mathbf{D}\mathbf{S}^{t}-\mathbf{Y}^{t})^{T}\mathbf{D}(\mathbb{B}_{S}\mathbf{Z}^{t}-\mathbf{S}^{t}))+\mu(\left\Vert \mathbb{B}_{S}\mathbf{Z}^{t}\right\Vert _{1}-\left\Vert \mathbf{S}^{t}\right\Vert _{1})}{\left\Vert \mathbf{D}(\mathbb{B}_{Q}\mathbf{Z}^{t}-\mathbf{S}^{t})\right\Vert _{F}^{2}}\right]_{0}^{1}.\label{eq:network-anomaly-S-stepsize}
\end{equation}
\setcounter{equation}{\value{MYtempeqncnt1}} \hrulefill  \end{figure*}

Consider the problem of estimating a low rank matrix $\mathbf{X}\in\mathbb{R}^{N\times K}$
and a sparse matrix $\mathbf{S}\in\mathbb{R}^{I\times K}$ from the
noisy measurement $\mathbf{Y}\in\mathbb{R}^{N\times K}$ which is
the output of a linear system:
\[
\mathbf{Y}=\mathbf{X}+\mathbf{DS}+\mathbf{V},
\]
where $\mathbf{D}\in\mathbb{R}^{N\times I}$ is known and $\mathbf{V}^{N\times K}$
is the unknown noise. The rank of $\mathbf{X}$ is much smaller than
$N$ and $K$, i.e, $\textrm{rank}(\mathbf{X})\ll\min(N,K)$, and
the support size of $\mathbf{S}$ is much smaller than $IK$, i.e.,
$\left\Vert \mathbf{S}\right\Vert _{0}\ll IK$. A commonly used regularization
function to promote the low rank is the nuclear norm $\left\Vert \mathbf{X}\right\Vert _{*}$,
but its has a cubic complexity which becomes unaffordable when the
problem dimension increases. It follows from the identity \cite{Burer2003,Recht2010}
$\left\Vert \mathbf{X}\right\Vert _{*}=\min_{(\mathbf{P},\mathbf{Q}):\mathbf{PQ=X}}\frac{1}{2}\left(\left\Vert \mathbf{P}\right\Vert _{F}^{2}+\left\Vert \mathbf{Q}\right\Vert _{F}^{2}\right)$
that the low rank matrix $\mathbf{X}$ can be written as the product
of two low rank matrices $\mathbf{P}\in\mathbb{R}^{N\times\rho}$
and $\mathbf{Q}\in\mathbb{R}^{\rho\times K}$ for a $\rho$ that is
usually much smaller than $N$ and $K$.

A natural measure for the estimation error is the least square loss
function augmented by regularization functions to promote the rank
sparsity of $\mathbf{X}$ and support sparsity of $\mathbf{S}$:
\begin{equation}
\underset{\mathbf{P},\mathbf{Q},\mathbf{S}}{\textrm{minimize}}\;\frac{1}{2}\left\Vert \mathbf{P}\mathbf{Q}\negthickspace+\negthickspace\mathbf{D}\mathbf{S}\negthickspace-\negthickspace\mathbf{Y}\right\Vert _{F}^{2}+\frac{\lambda}{2}\left(\left\Vert \mathbf{P}\right\Vert _{F}^{2}\negthickspace+\negthickspace\left\Vert \mathbf{Q}\right\Vert _{F}^{2}\right)+\mu\left\Vert \mathbf{S}\right\Vert _{1},\label{eq:eq:rank-problem-formulation}
\end{equation}
where the matrix factorization $\mathbf{X}=\mathbf{PQ}$ has been
used and it does not incur any estimation error under some sufficient
conditions specified in \cite[Prop. 1]{Mardani2013}. This nonconvex
optimization problem is a special case of (\ref{eq:problem-formulation})
obtained by setting $f(\mathbf{P},\mathbf{Q},\mathbf{S})\triangleq\frac{1}{2}\left\Vert \mathbf{P}\mathbf{Q}+\mathbf{D}\mathbf{S}-\mathbf{Y}\right\Vert _{F}^{2}+\frac{\lambda}{2}\left(\left\Vert \mathbf{P}\right\Vert _{F}^{2}+\left\Vert \mathbf{Q}\right\Vert _{F}^{2}\right)$
and $g(\mathbf{S})\triangleq\mu\left\Vert \mathbf{S}\right\Vert _{1}$.
Note that $\nabla f$ is not Lipschitz continuous. To see this, consider
the scalar case: its gradient $\nabla_{P}f=(PQ+DS-Y)Q$ and $\left|(P^{'}Q^{2}+DS-Y)-(P^{''}Q^{2}+DS-Y)\right|\leq Q^{2}\left|P'-P''\right|$
while the unconstrained $Q$ can be unbounded (it is however block
Lipschitz continuous).

The problem formulation (\ref{eq:eq:rank-problem-formulation}) plays
an important role in the network anomaly detection problem in \cite{Mardani2013b}.
A parallel SCA algorithm in the essence of Alg. \ref{alg:Successive-approximation-method}
was proposed in \cite{Yang_NonconvexRegularization,Yang_Rank_ICASSP2018}
to solve problem (\ref{eq:eq:rank-problem-formulation}), where $\mathbf{P}$,
$\mathbf{Q}$ and $\mathbf{S}$ are updated simultaneously at each
iteration. However, it assumes the memory capacity is large enough
to store the whole data set and all intermediate variables generated
at each iteration.

In this section, we apply the BSCA algorithm proposed in Sec. \ref{sec:Block-SCA}
to solve problem (\ref{eq:eq:rank-problem-formulation}). Define $\mathbf{Z}\triangleq(\mathbf{P},\mathbf{Q},\mathbf{S})$
and assume for simplicity the cyclic update rule. As $f(\mathbf{P},\mathbf{Q},\mathbf{S})$
is individually convex in $\mathbf{P}$, $\mathbf{Q}$ and $\mathbf{S}$,
the approximation function w.r.t. one block variable is obtained by
fixing other block variables (cf. (\ref{eq:approximation-best-response-1})-(\ref{eq:approximation-best-response-2})):
\begin{subequations}\label{eq:network-anomaly-approximate-function}
\begin{align}
\widetilde{f}(\mathbf{P};\mathbf{Z}^{t}) & =f(\mathbf{P},\mathbf{Q}^{t},\mathbf{S}^{t}),\label{eq:network-anomaly-approximate-function-P}\\
\widetilde{f}(\mathbf{Q};\mathbf{Z}^{t}) & =f(\mathbf{P}^{t},\mathbf{Q},\mathbf{S}^{t}),\label{eq:network-anomaly-approximate-function-Q}\\
\widetilde{f}(\mathbf{S};\mathbf{Z}^{t}) & =\sum_{i,k}f(\mathbf{P}^{t},\mathbf{Q}^{t},s_{i,k},(s_{j,k}^{t})_{j\neq i},(\mathbf{s}_{j}^{t})_{j\neq k})\nonumber \\
 & =\sum_{i,k}\frac{1}{2}\left\Vert \mathbf{P}^{t}\mathbf{q}_{k}^{t}+\mathbf{d}_{i}s_{i,k}+{\textstyle \sum_{j\neq i}}\mathbf{d}_{j}s_{j,k}^{t}-\mathbf{y}_{k}\right\Vert _{2}^{2},\label{eq:network-anomaly-approximate-function-S}
\end{align}
\end{subequations}where $\mathbf{q}_{k}$, $\mathbf{d}_{k}$ and
$\mathbf{y}_{k}$ in (\ref{eq:network-anomaly-approximate-function-S})
is the $k$-th column of $\mathbf{Q}$, $\mathbf{D}$ and $\mathbf{Y}$,
respectively. Note that it may be tempting to approximate $f(\mathbf{P},\mathbf{Q},\mathbf{S})$
w.r.t. $\mathbf{S}$ by $f(\mathbf{P}^{t},\mathbf{Q}^{t},\mathbf{S})$,
but the resulting approximation subproblem does not have a closed-form
solution and must be solved by iterative algorithms.

If the block variable $\mathbf{P}$ or $\mathbf{Q}$ is updated at
iteration $t$, the approximation subproblem is\begin{subequations}\label{eq:network-anomaly-approximate-problem}
\begin{align}
\mathbb{B}_{P}\mathbf{Z}^{t} & =\underset{\mathbf{P}}{\arg\min}\;\widetilde{f}(\mathbf{P};\mathbf{Z}^{t})\nonumber \\
 & =(\mathbf{Y}-\mathbf{D}\mathbf{S}^{t})(\mathbf{Q}^{t})^{T}(\mathbf{Q}^{t}(\mathbf{Q}^{t})^{T}+\lambda\mathbf{I})^{-1}.\label{eq:network-anomaly-approximate-problem-P}
\end{align}
or
\begin{align}
\mathbb{B}_{Q}\mathbf{Z}^{t} & =\underset{\mathbf{Q}}{\arg\min}\;\widetilde{f}(\mathbf{Q};\mathbf{Z}^{t})\nonumber \\
 & =((\mathbf{P}^{t})^{T}\mathbf{P}^{t}+\lambda\mathbf{I})^{-1}(\mathbf{P}^{t})^{T}(\mathbf{Y}-\mathbf{D}\mathbf{S}^{t}),\label{eq:network-anomaly-approximate-problem-Q}
\end{align}
respectively. When the block variable $\mathbf{S}$ is updated, the
approximation subproblem is
\begin{align}
\mathbb{B}_{S}\mathbf{Z}^{t} & =\underset{\mathbf{S}}{\arg\min}\;\widetilde{f}(\mathbf{S};\mathbf{Z}^{t})+g(\mathbf{S})\nonumber \\
 & =\mathbf{d}(\mathbf{D}^{T}\mathbf{D})^{-1}\cdot\nonumber \\
 & \qquad\mathcal{S}_{\mu}\left(\mathbf{d}(\mathbf{D}^{T}\mathbf{D})\mathbf{S}^{t,\tau}-\mathbf{D}^{T}(\mathbf{D}\mathbf{S}^{t,\tau}-\mathbf{Y}^{t}+\mathbf{P}^{t}\mathbf{Q}^{t})\right),\label{eq:network-anomaly-approximate-problem-S}
\end{align}
\end{subequations}where $\mathcal{S}_{a}[\mathbf{X}]\triangleq\max(\mathbf{X}-a\mathbf{I},\mathbf{0})-\max(-\mathbf{X}-a\mathbf{I},\mathbf{0})$
is the soft-thresholding operator. The next point $\mathbf{S}^{t+1}$
is defined as
\begin{equation}
\mathbf{S}^{t+1}=\mathbf{S}^{t}+\gamma^{t}(\mathbb{B}_{S}\mathbf{Z}^{t}-\mathbf{S}^{t}).\label{eq:network-anomaly-S-update}
\end{equation}
The stepsize $\gamma^{t}$ can be obtained by performing the exact
line search along the coordinate of $\mathbf{S}$ over $f(\mathbf{P}^{t},\mathbf{Q}^{t},\mathbf{S})+g(\mathbf{S})$:
\[
\gamma^{t}=\underset{0\leq\gamma\leq1}{\arg\min}\left\{ \begin{array}{l}
f(\mathbf{P}^{t},\mathbf{Q}^{t},\mathbf{S}^{t}+\gamma(\mathbb{B}_{S}\mathbf{Z}^{t}-\mathbf{S}^{t}))\smallskip\\
+\gamma(g(\mathbb{B}_{S}\mathbf{Z}^{t})-g(\mathbf{S}^{t}))
\end{array}\right\} .
\]
It has a closed-form expression given at the top of this page.\addtocounter{equation}{1}

The above steps are summarized in Alg. \ref{alg:network-anomaly-detection}.
Note that when updating $\mathbf{P}$ or $\mathbf{Q}$, we have used
a constant unit stepsize because the approximation function $\widetilde{f}(\mathbf{P};\mathbf{Z}^{t})$
in (\ref{eq:network-anomaly-approximate-function-P}) and $\widetilde{f}(\mathbf{Q};\mathbf{Z}^{t})$
in (\ref{eq:network-anomaly-approximate-function-Q}) is a (trivial)
global upper bound of $f(\mathbf{P},\mathbf{Q}^{t},\mathbf{S}^{t})$
and $f(\mathbf{P}^{t},\mathbf{Q},\mathbf{S}^{t})$, respectively (see
the discussion for (\ref{eq:variable-update-upper-bound})). It follows
from Theorem \ref{thm:block-SCA-convergence} that every limit point
of the sequence $\{\mathbf{P}^{t},\mathbf{Q}^{t},\mathbf{S}^{t}\}$
generated by Alg. \ref{alg:network-anomaly-detection} is a stationary
point of (\ref{eq:eq:rank-problem-formulation}).

\begin{algorithm}[t]
\textbf{Initialization: }$t=0$, $\mathbf{Z}^{0}=(\mathbf{P}^{0},\mathbf{Q}^{0},\mathbf{S}^{0})$
(arbitrary but fixed).

Repeat the following steps until convergence:

\begin{enumerate}

\item[\textbf{S1: }] Choose a block variable ($\mathbf{P}$, $\mathbf{Q}$
or $\mathbf{S}$) according to either the cyclic update rule or the
random update rule.

\end{enumerate}

If $\mathbf{P}$ is selected:

\begin{enumerate}

\item[\textbf{S2: }] $\mathbf{Q}^{t+1}=\mathbf{Q}^{t}$, $\mathbf{S}^{t+1}=\mathbf{S}^{t}$,
and $\mathbf{P}^{t+1}=\mathbb{B}_{P}\mathbf{Z}^{t}$ defined in (\ref{eq:network-anomaly-approximate-problem-P}).

\end{enumerate}

If $\mathbf{Q}$ is selected:

\begin{enumerate}

\item[\textbf{S2: }] $\mathbf{P}^{t+1}=\mathbf{P}^{t}$, $\mathbf{S}^{t+1}=\mathbf{S}^{t}$
and $\mathbf{Q}^{t+1}=\mathbb{B}_{Q}\mathbf{Z}^{t}$ defined in (\ref{eq:network-anomaly-approximate-problem-Q}).

\end{enumerate}

If $\mathbf{S}$ is selected:

\begin{enumerate}

\item[\textbf{S2: }] $\mathbf{P}^{t+1}=\mathbf{P}^{t}$, $\mathbf{Q}^{t+1}=\mathbf{Q}^{t}$,
and $\mathbf{S}^{t+1}$ is obtained by the following steps:

\textbf{}%
\colorbox{lightgray}{\begin{minipage}[t]{0.88\columnwidth}%
\begin{enumerate}

\item[\textbf{S2.1:}] Compute $\mathbb{B}_{S}\mathbf{Z}^{t}$ according
to (\ref{eq:network-anomaly-approximate-problem-S}).

\item[\textbf{S2.2:}] Compute the stepsize $\gamma^{t}$ by the exact
line search (\ref{eq:network-anomaly-S-stepsize}).

\item[\textbf{S2.3:}] Update $\mathbf{S}^{t+1}$ according to (\ref{eq:network-anomaly-S-update}).

\end{enumerate}%
\end{minipage}}

\end{enumerate}

\begin{enumerate}

\item[\textbf{S3: }] $t\leftarrow t+1$ and go to \textbf{S1}.

\end{enumerate}

\caption{\label{alg:network-anomaly-detection}The exact block successive convex
approximation algorithm for network anomaly detection problem (\ref{eq:eq:rank-problem-formulation})}
\end{algorithm}

We remark that Alg. \ref{alg:network-anomaly-detection} enjoys i)
low complexity as all variable updates can be performed by closed-form
expressions; ii) easy implementation as the three block variables
are updated sequentially and only a single processor is needed; and
iii) fast convergence as when a particular block variable is updated,
the most recent updates of previous blocks are exploited. Although
the gradient of the objective function w.r.t. each block variable
is Lipschitz continuous, the proposed algorithm does not have any
hyperparameters that are dependent on the typically unknown Lipschitz
continuity constant.

\begin{figure*}[t]
\center

\center\includegraphics[bb=100bp 20bp 950bp 480bp,clip,scale=0.6]{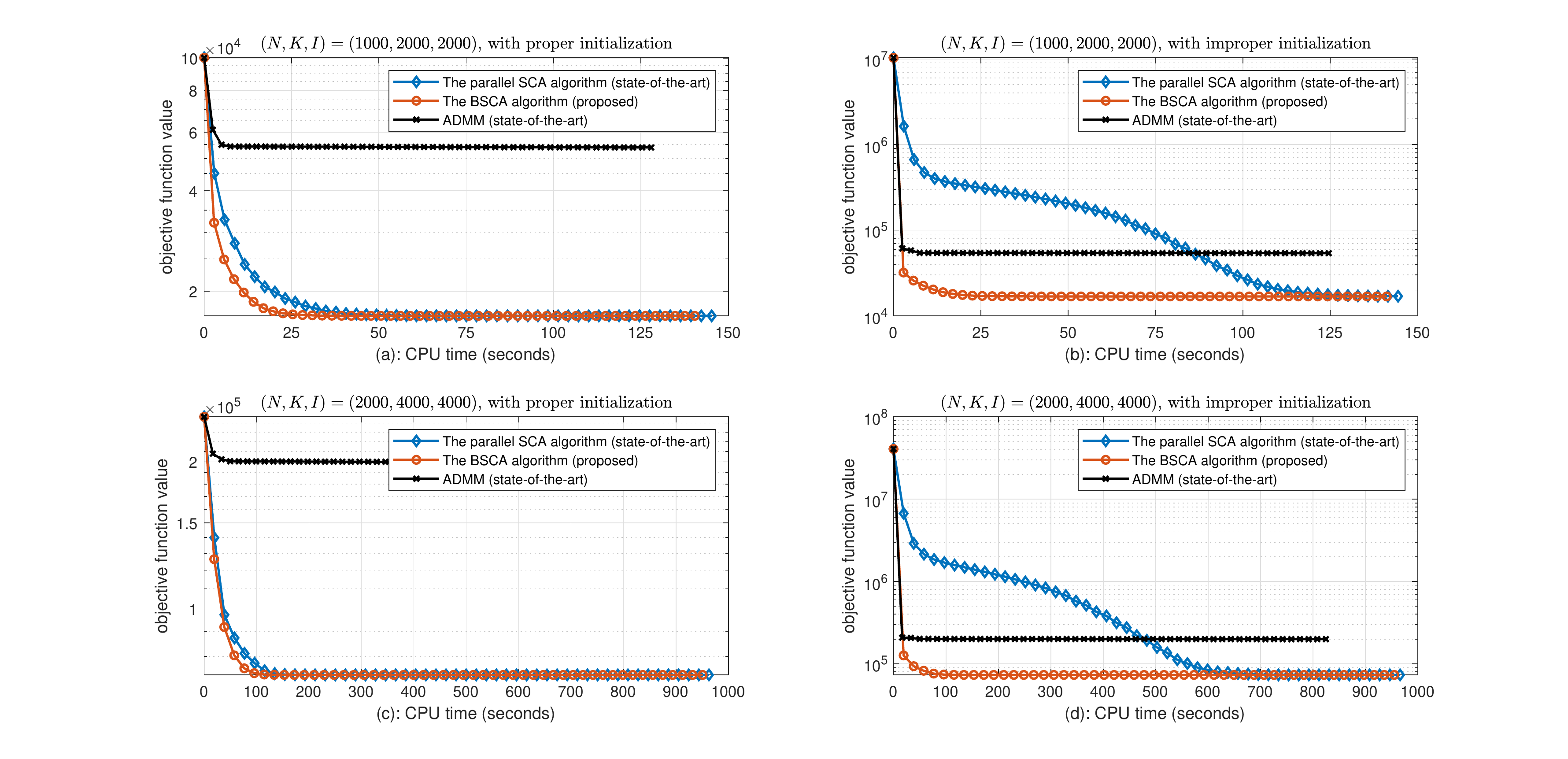}

\caption{\label{fig:Rank-Minimization}Joint estimation of low-rank and sparse
signals: The achieved objective function value versus the CPU time.}
\end{figure*}

\textbf{Simulations. }All simulations in this paper are carried out
under Matlab R2019a on a laptop equipped with a Windows 7 64-bit operating
system, an Intel i5-3210 2.50GHz CPU with 4 logical processors, and
a 8GB RAM. Although the proposed updates involve linear algebraic
operations only, we do not write low-level program to directly call
the processors and parallelize the proposed algorithms. Instead, we
rely on the computer compiler and numerical libraries (for example
LAPACK), both of which are nowadays highly optimized and well integrated
for parallel computations in computing softwares such as Matlab and
coding languages such as Python, to parallelize the linear algebraic
operations.

The simulation parameters are set as follows. $(N,K,I)=(1000,2000,2000)$
or $(2000,4000,4000)$, $\rho=5$. The regularization parameters $\lambda=0.25\left\Vert \mathbf{Y}\right\Vert $
and $\mu=2\cdot10^{-4}\left\Vert \mathbf{D}^{T}\mathbf{Y}\right\Vert _{\infty}$.
The elements of $\mathbf{D}$ are first generated according to the
normal distribution, and each row is then normalized to unity. The
elements of $\mathbf{V}$ follow the Gaussian distribution with mean
0 and variance $10^{-4}$. The density of $\mathbf{S}$ is 0.05 and
its nonzero elements are generated according to the normal distribution.
We set $\mathbf{Y}=\mathbf{PQ}+\mathbf{DS}+\mathbf{V}$, where $\mathbf{P}$
and \textbf{$\mathbf{Q}$} are generated randomly following the Gaussian
distribution $\mathcal{N}(0,100/I)$ and $\mathcal{N}(0,100/K)$,
respectively. The simulation results are averaged over 20 realizations.

In Fig. \ref{fig:Rank-Minimization} the achieved objective function
value versus the CPU time of the parallel SCA algorithm \cite{Yang_NonconvexRegularization},
ADMM algorithm \cite{Mardani2013} and the proposed BSCA algorithm
is plotted. The marker on the curve represents an iteration. All algorithms
start with two different initializations: ``proper initialization''
if $\mathbf{P}^{0}$ and $\mathbf{Q}^{0}$ are generated in the same
way as the real $\mathbf{P}$ and $\mathbf{Q}$, i.e., all elements
follow the Gaussian distribution $\mathcal{N}(0,100/I)$ and $\mathcal{N}(0,100/K)$,
or ``improper initialization'' if they are generated randomly following
the standard normal distribution $\mathcal{N}(0,1)$ and $\mathcal{N}(0,1)$.
For the parallel SCA algorithm, the code is divided into blocks and
the parallelizable blocks are executed sequentially.

From Fig. \ref{fig:Rank-Minimization} we can draw several observations.

Both BSCA and SCA algorithms converge to the same objective function
value, which is notably better than the value to which the ADMM algorithm
converges to. Although the ADMM appears to be convergent in the simulations,
it does not have a guaranteed convergence.

We see from Fig. \ref{fig:Rank-Minimization} that the BSCA algorithm
exhibits a faster convergence in terms of the CPU time than naively
dividing the parallel SCA algorithm into blocks and executing the
parallel blocks sequentially, especially when the initial point is
far away from the optimal point. This consolidates the intuition that
exploiting the most recent update of previous block variables is beneficial
and could significantly accelerate the convergence.

Comparing Fig. \ref{fig:Rank-Minimization} (a) with Fig. \ref{fig:Rank-Minimization}
(b) and Fig. \ref{fig:Rank-Minimization} (c) with Fig. \ref{fig:Rank-Minimization}
(d), we see that the SCA algorithm is more sensitive to the choice
of the initial point. By contrast, the BSCA algorithm converges to
the optimal point in the same number of iterations (which can be counted
by the number of markers) and the same CPU time.

When increasing the problem dimension from $(N,K,I)=(1000,2000,2000)$
in Fig. \ref{fig:Rank-Minimization} (a)-(b) to $(N,K,I)=(2000,4000,4000)$
in Fig. \ref{fig:Rank-Minimization} (c)-(d), we see that the BSCA
algorithm still converges to an accurate solution within 10 iterations
(the CPU time increases as the higher problem dimension leads to higher
computational complexity per iteration). Therefore the BSCA algorithm
scales very well.

\subsection{Quadratic Inverse Problems and Phase Retrieval}

In phase retrieval problems, we are given a number of magnitude measurements
that are of the following form
\[
y_{n}\approx(\mathbf{a}_{n}^{T}\mathbf{x}_{0})^{2},n=1,\ldots,N,
\]
where $\mathbf{x}_{0}$ is the unknown sparse signal, $\mathbf{a}_{n}$
is a known sampling vector\footnote{For simplicity we assume $\mathbf{x}_{0}$ and $\mathbf{a}_{n}$ are
real-valued, but all results can be generalized to the complex-valued
case.}, and $N$ is the number of observations. To estimate $\mathbf{x}_{0}$
from the noisy magnitude measurements $(y_{n})_{n=1}^{N}$, one of
the most popular approaches is optimization-based approach, which
amounts to solving a nonconvex quadratic inverse problem
\begin{align}
\underset{\mathbf{x}\in\mathbb{R}^{I}}{\textrm{minimize}}\; & \frac{1}{4}\sum_{n=1}^{N}\left((\mathbf{a}_{n}^{T}\mathbf{x})^{2}-y_{n}\right)^{2}+\mu\left\Vert \mathbf{x}\right\Vert _{1}.\label{eq:QI-formulation}
\end{align}
Quadratic inverse problems are also referred to as the phase retrieval
problem \cite{Candes2015} and it is an instance of (\ref{eq:problem-formulation})
with the decomposition
\[
f(\mathbf{x})=\frac{1}{4}\sum_{n=1}^{N}\left((\mathbf{a}_{n}^{T}\mathbf{x})^{2}-y_{n}\right)^{2},\textrm{ and }g(\mathbf{x})=\mu\left\Vert \mathbf{x}\right\Vert _{1}.
\]
Note that $\nabla f$ is not block Lipschitz continuous. To see this,
consider the special case $f(x)=\frac{1}{4}(x^{2}-y)^{2}$: its gradient
is $x(x^{2}-y)$ and thus not Lipschitz continuous.

Define $l_{n}(\mathbf{x})\triangleq(\mathbf{a}_{n}^{T}\mathbf{x})^{2}-y_{n}$
and rewrite $f$ as the composition of functions $f(\mathbf{x})=\sum_{n=1}^{N}\frac{1}{4}(l_{n}(\mathbf{x}))^{2}$.
To apply the proposed BSCA algorithm, we first approximate $f(\mathbf{x})$
by the partial linearization approximation (see (\ref{eq:approximation-partial-linearization})),
that is,
\begin{align}
\widetilde{f}(\mathbf{x}_{k};\mathbf{x}^{t})=\; & \frac{1}{4}\sum_{n=1}^{N}(l_{n}(\mathbf{x}^{t})+(\mathbf{x}_{k}-\mathbf{x}_{k}^{t})^{T}\nabla_{k}l_{n}(\mathbf{x}^{t}))^{2}\nonumber \\
 & +\frac{c_{k}^{t}}{2}\left\Vert \mathbf{x}_{k}-\mathbf{x}_{k}^{t}\right\Vert ^{2},\nonumber \\
=\; & \frac{1}{2}\mathbf{x}_{k}^{T}\mathbf{D}_{k}^{t}\mathbf{x}_{k}-\mathbf{x}_{k}^{T}\mathbf{b}_{k}^{t},\label{eq:QI-outer-approximation-function}
\end{align}
where $c_{k}^{t}$ is a positive scalar, and
\begin{align*}
\mathbf{D}_{k}^{t} & \triangleq2\mathbf{A}_{k}\textrm{diag}(\mathbf{A}^{T}\mathbf{x}^{t})\textrm{diag}(\mathbf{A}^{T}\mathbf{x}^{t})\mathbf{A}_{k}^{T}+c_{k}^{t}\mathbf{I},\\
\mathbf{b}_{k}^{t} & \triangleq\mathbf{D}_{k}^{t}\mathbf{x}_{k}^{t}-\mathbf{A}_{k}((\mathbf{A}^{T}\mathbf{x}^{t})\circ\mathbf{l}(\mathbf{x}^{t})),
\end{align*}
with $\mathbf{A}_{k}\in\mathbb{R}^{I_{k}\times N}$, $\sum_{k=1}^{K}I_{k}=I$
and
\begin{align*}
\mathbf{A} & \triangleq\left[\begin{array}{ccccc}
\mathbf{a}_{1} & \ldots & \mathbf{a}_{n} & \ldots & \mathbf{a}_{N}\end{array}\right]=\left[\begin{array}{c}
\mathbf{A}_{1}\\
\vdots\\
\mathbf{A}_{k}\\
\vdots\\
\mathbf{A}_{K}
\end{array}\right]\in\mathbb{R}^{I\times N}.
\end{align*}
It can be verified that
\[
\nabla_{k}\widetilde{f}(\mathbf{x}^{t};\mathbf{x}^{t})=\nabla_{k}f(\mathbf{x}^{t})=\mathbf{A}_{k}((\mathbf{A}^{T}\mathbf{x}^{t})\circ\mathbf{l}(\mathbf{x}^{t})).
\]

The (outer-layer) approximation subproblem is
\begin{align}
\underset{\mathbf{x}}{\textrm{minimize}}\quad & \widetilde{f}(\mathbf{x};\mathbf{x}^{t})+g(\mathbf{x}).\label{eq:QI-approx-outer}
\end{align}
This problem however does not have a closed-form solution and we solve
it inexactly by running the SCA algorithm for a finite number of iterations
in the inner layer. For the inner-layer approximation function, as
$\widetilde{f}(\mathbf{x}_{k};\mathbf{x}^{t})$ is strictly convex
in $\mathbf{x}_{k}$, we adopt the best-response approximation: given
$\mathbf{x}_{k}^{t,\tau}$ at iteration $\tau$ of the inner layer,
\[
\widetilde{f}_{k}^{i}(\mathbf{x}_{k};\mathbf{x}_{k}^{t,\tau},\mathbf{x}^{t})=\sum_{i_{k}=1}^{I_{k}}\widetilde{f}(x_{i_{k}},(x_{j_{k}}^{t,\tau})_{j_{k}\neq i_{k}};\mathbf{x}^{t}).
\]
The inner-layer approximation subproblem has a closed-form solution
\begin{align}
\mathbb{B}_{k}\overline{\mathbf{x}}^{t,\tau} & =\underset{\mathbf{x}_{k}}{\arg\min}\;\widetilde{f}_{k}^{i}(\mathbf{x}_{k};\mathbf{x}_{k}^{t,\tau},\mathbf{x}^{t})+g_{k}(\mathbf{x}_{k})\nonumber \\
 & =S_{\mu\mathbf{d}(\mathbf{D}_{k}^{t})^{-1}}\left(\mathbf{x}^{t,\tau}-\frac{\mathbf{D}_{k}^{t}\mathbf{x}^{t,\tau}-\mathbf{b}_{k}^{t}}{\mathbf{d}(\mathbf{D}_{k}^{t})}\right),\label{eq:QI-approximation-subproblem-inner}
\end{align}
where $S$ is the soft-thresholding operator, and the vector division
is understood to be an element-wise operation.

\begin{algorithm}[t]
\textbf{Initialization: }$t=0$, $\mathbf{x}^{0}$ (nonzero, arbitrary
but fixed).

Repeat the following steps until convergence:

\textbf{S1: }Select the block variable $\mathbf{x}_{k}$ to be updated
according to either the cyclic update rule or the random update rule.

\textbf{S2: }Compute $\widetilde{\mathbf{x}}_{k}^{t}$ by the following
steps:

\hspace{0.5cm}\textbf{}%
\colorbox{lightgray}{\begin{minipage}[t]{0.88\columnwidth}%
\begin{enumerate}

\item[\textbf{S1.0:}] Set $\tau=0$ and $\mathbf{x}_{k}^{t,0}=\mathbf{x}_{k}^{t}$.

\item[\textbf{S1.1:}] Compute $\mathbb{B}_{k}\overline{\mathbf{x}}^{t,\tau}$
according to (\ref{eq:QI-approximation-subproblem-inner}).

\item[\textbf{S1.2:}] Compute the stepsize $\gamma^{t,\tau}$ according
to (\ref{eq:QI-stepsize-inner}).

\item[\textbf{S1.3:}] Update $\mathbf{x}_{k}^{t,\tau+1}$ according
to (\ref{eq:QI-update-inner}).

\item[\textbf{S1.4:}] If $\tau+1=\bar{\tau}_{t}$, $\widetilde{\mathbf{x}}_{k}^{t}=\mathbf{x}_{k}^{t,\bar{\tau}_{t}}$
and go to \textbf{S3}. Otherwise $\tau\leftarrow\tau+1$ and go to
\textbf{S2.1}.

\end{enumerate}%
\end{minipage}}

\textbf{S3: }Compute the stepsize $\gamma^{t}$ by the exact line
search (\ref{eq:QI-stepsize-outer}).

\textbf{S4: }Update $\mathbf{x}^{t+1}$ according to (\ref{eq:QI-update-outer}).

\textbf{S5: }$t\leftarrow t+1$ and go to \textbf{S1}.

\caption{\label{alg:Quadratic-Inverse}The inexact block successive convex
approximation algorithm for quadratic inverse problem (\ref{eq:QI-formulation})}
\end{algorithm}

Given the descent direction $\mathbb{B}_{k}\overline{\mathbf{x}}^{t,\tau}-\mathbf{x}_{k}^{t,\tau}$,
we refine $\mathbf{x}_{k}^{t,\tau}$ as
\begin{equation}
\mathbf{x}_{k}^{t,\tau+1}=\mathbf{x}_{k}^{t,\tau}+\gamma^{t,\tau}(\mathbb{B}_{k}\overline{\mathbf{x}}^{t,\tau}-\mathbf{x}_{k}^{t,\tau}),\label{eq:QI-update-inner}
\end{equation}
and the stepsize is obtained by performing the exact line search,
which has a simple analytical expression
\begin{align}
\gamma^{t,\tau} & =\underset{0\leq\gamma\leq1}{\arg\min}\left\{ \begin{array}{l}
\widetilde{f}(\mathbf{x}_{k}^{t}+\gamma(\mathbb{B}_{k}\overline{\mathbf{x}}^{t,\tau}-\mathbf{x}_{k}^{t,\tau});\mathbf{x}^{t})\smallskip\\
g(\mathbf{x}_{k}^{t,\tau})+\gamma(g(\mathbb{B}_{k}\overline{\mathbf{x}}^{t,\tau})-g(\mathbf{x}_{k}^{t,\tau}))
\end{array}\right\} \nonumber \\
= & \left[-\frac{(\mathbf{D}_{k}^{t}\mathbf{x}_{k}^{t,\tau}-\mathbf{b}_{k}^{t})\triangle\mathbf{x}_{k}^{t,\tau}+\mu\left(\left\Vert \mathbb{B}_{k}\overline{\mathbf{x}}^{t,\tau}\right\Vert _{1}-\left\Vert \mathbf{x}_{k}^{t,\tau}\right\Vert _{1}\right)}{(\triangle\mathbf{x}_{k}^{t,\tau})^{T}\mathbf{D}_{k}^{t}\triangle\mathbf{x}_{k}^{t,\tau}}\right]_{0}^{1},\label{eq:QI-stepsize-inner}
\end{align}
with $\triangle\mathbf{x}_{k}^{t,\tau}\triangleq\mathbb{B}_{k}\overline{\mathbf{x}}^{t,\tau}-\mathbf{x}_{k}^{t,\tau}$.
After repeating the above steps for a finite number of iterations,
we obtain an inexact solution of the outer-layer approximation subproblem
(\ref{eq:QI-approx-outer}), which we denote as $\widetilde{\mathbf{x}}_{k}^{t}$.

\begin{figure*}[tbh]
\center

\includegraphics[bb=75bp 5bp 820bp 330bp,clip,scale=0.68]{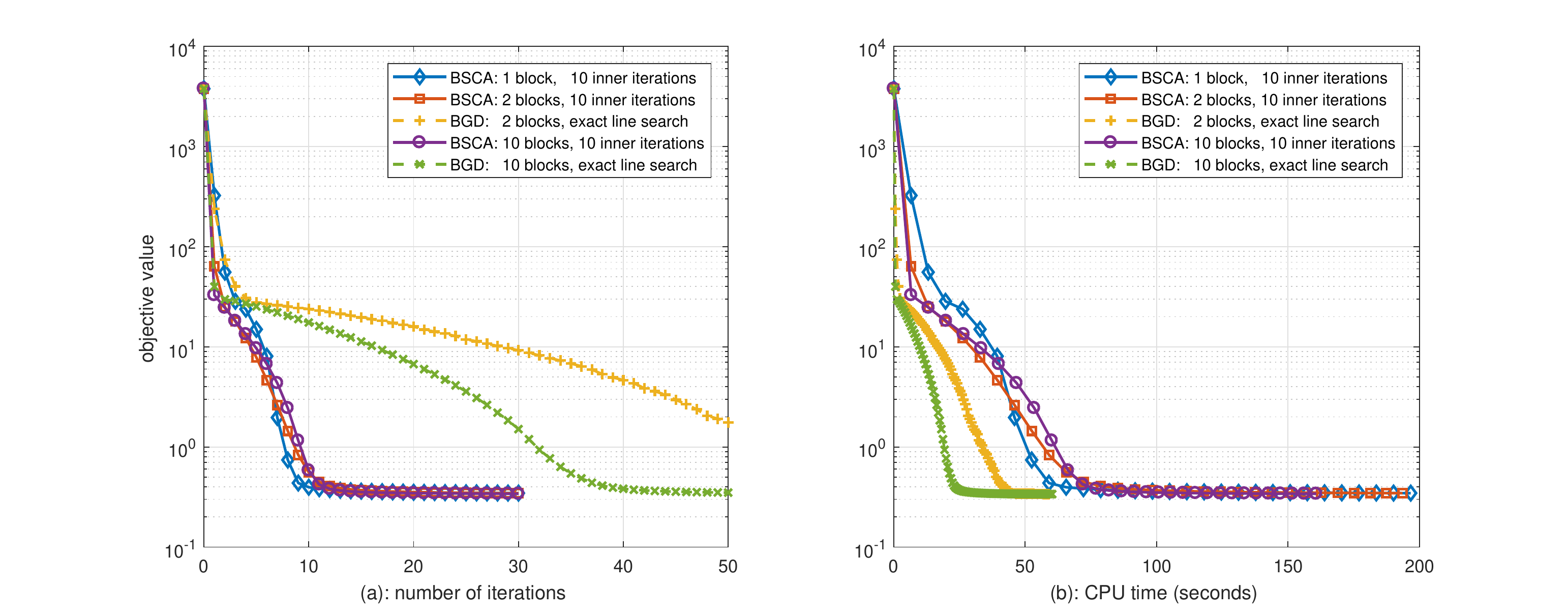}

\includegraphics[bb=75bp 5bp 820bp 330bp,clip,scale=0.68]{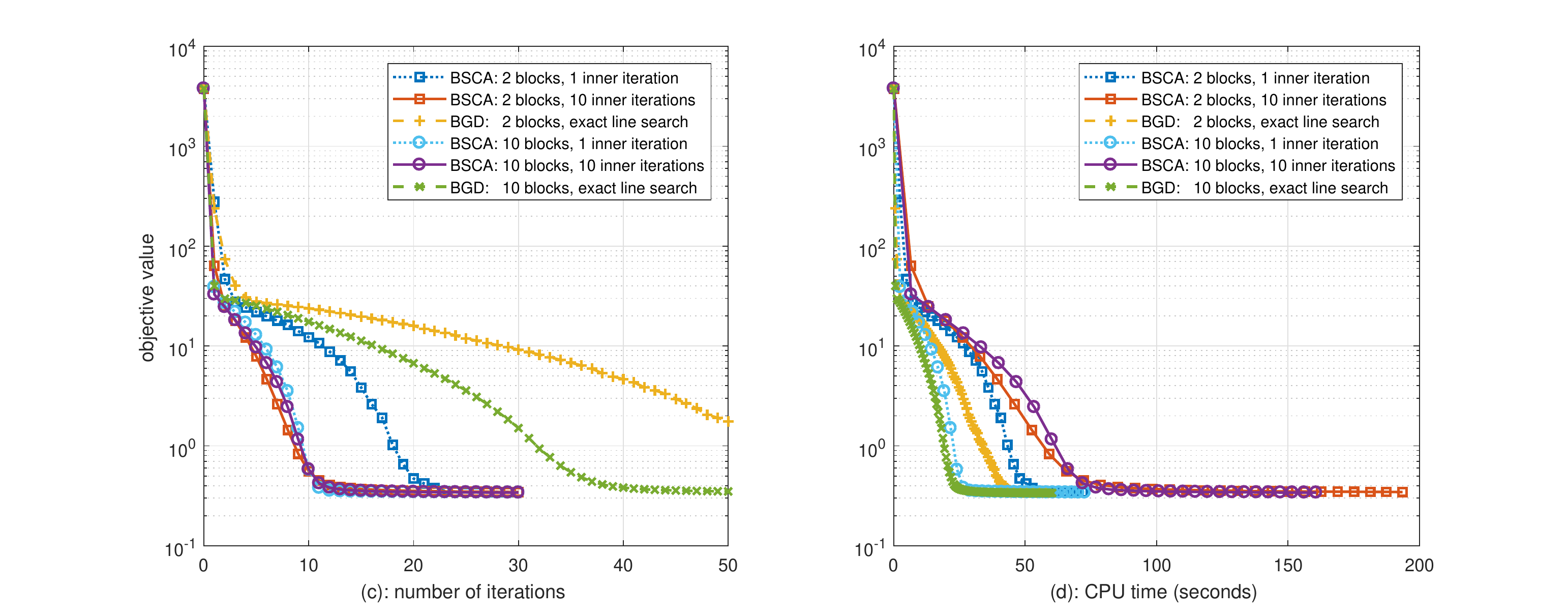}\caption{\label{fig:QI}Phase retrieval: objective value versus the number
of iterations and the CPU time}
\end{figure*}

Since $\widetilde{\mathbf{x}}_{k}^{t}-\mathbf{x}_{k}^{t}$ is a descent
direction of $f+g$ along the coordinate of $\mathbf{x}_{k}$, we
are ready to refine $\mathbf{x}^{t}$:
\begin{equation}
\mathbf{x}_{k}^{t+1}=\mathbf{x}_{k}^{t}+\gamma^{t}(\widetilde{\mathbf{x}}_{k}^{t}-\mathbf{x}_{k}^{t}),\label{eq:QI-update-outer}
\end{equation}
and $\mathbf{x}_{j}^{t+1}=\mathbf{x}_{j}^{t}$ for all $j\neq k$.
We choose to compute the stepsize $\gamma^{t}$ in the outer layer
by the exact line search
\begin{align}
\gamma^{t} & =\underset{0\leq\gamma\leq1}{\arg\min}\left\{ \begin{array}{l}
f(\mathbf{x}_{k}^{t}+\gamma(\widetilde{\mathbf{x}}_{k}^{t}-\mathbf{x}_{k}^{t}),\mathbf{x}_{-k}^{t})\smallskip\\
g_{k}(\mathbf{x}_{k}^{t})+\gamma(g_{k}(\widetilde{\mathbf{x}}_{k}^{t})-g_{k}(\mathbf{x}_{k}^{t}))
\end{array}\right\} \nonumber \\
 & =\underset{0\leq\gamma\leq1}{\arg\min}\left\{ \frac{1}{4}v_{4}\gamma^{4}+\frac{1}{3}v_{3}\gamma^{3}+\frac{1}{2}v_{2}\gamma^{2}+v_{1}\gamma\right\} ,\label{eq:QI-stepsize-outer}
\end{align}
where
\begin{align*}
v_{4}=\; & \left\Vert (\mathbf{A}_{k}^{T}\mathbf{x}_{k}^{t})^{2}\right\Vert _{2}^{2},\\
v_{3}=\; & 3(\mathbf{A}^{T}\mathbf{x}^{t})^{T}(\mathbf{A}_{k}^{T}\triangle\mathbf{x}_{k}^{t})^{3},\\
v_{2}=\; & (3(\mathbf{A}^{T}\mathbf{x}^{t})^{2}-\mathbf{y})^{T}(\mathbf{A}_{k}^{T}\triangle\mathbf{x}_{k}^{t})^{2},\\
v_{1}=\; & (\mathbf{A}_{k}^{T}\triangle\mathbf{x}_{k}^{t})^{T}\left((\mathbf{A}^{T}\mathbf{x}^{t})^{3}-(\mathbf{A}^{T}\mathbf{x}^{t})\circ\mathbf{y}\right)\\
 & +\mu(\left\Vert \widetilde{\mathbf{x}}_{k}^{t}\right\Vert _{1}-\left\Vert \mathbf{x}_{k}^{t}\right\Vert _{1}),
\end{align*}
with $\triangle\mathbf{x}_{k}^{t}\triangleq\widetilde{\mathbf{x}}_{k}^{t}-\mathbf{x}_{k}^{t}$.
Solving the optimization problem in (\ref{eq:QI-stepsize-outer})
is equivalent to finding the nonnegative real root of a third-order
polynomial. By the Cardano's method, $\gamma^{t}$ has an analytical
expression\begin{subequations}\label{eq:rank-line-search-proposed-closed-form}
\begin{align}
\gamma^{t} & =[\bar{\gamma}^{t}]_{0}^{1},\label{eq:rank-line-search-proposed-closed-form-1}\\
\bar{\gamma}^{t} & =\sqrt[3]{\Sigma_{1}+\sqrt{\Sigma_{1}^{2}+\Sigma_{2}^{3}}}+\sqrt[3]{\Sigma_{1}-\sqrt{\Sigma_{1}^{2}+\Sigma_{2}^{3}}}-\frac{v_{3}}{3v_{4}},\label{eq:rank-line-search-proposed-closed-form-2}
\end{align}
\end{subequations}where $\Sigma_{1}\triangleq-(v_{3}/3v_{4})^{3}+v_{3}v_{1}/6v_{4}^{2}-v_{1}/2v_{4}$
and $\Sigma_{2}\triangleq v_{1}/3v_{4}-(v_{3}/3v_{4})^{2}$. Note
that in (\ref{eq:rank-line-search-proposed-closed-form-2}), the right
hand side contains three values (two of them can attain complex numbers),
and the equal sign must be interpreted as assigning the smallest real
nonnegative values.

The above steps are summarized in Alg. \ref{alg:Quadratic-Inverse}
and it has several notable advantages. Firstly, it has a guaranteed
convergence to a stationary point, although the gradient of the smooth
function $f$ is not (block) Lipschitz continuous. Secondly, it exhibits
a fast convergence as the approximation function preserves the problem
structure to a large extent. Besides, it enables sequential block
update and is suitable for hardware with limited memory and/or processing
capability. Furthermore, it has low complexity as all updates have
analytical expressions.

\textbf{Simulations}. In our numerical simulations the dimension of
$\mathbf{A}$ is $5000\times20000$: all of its elements are generated
randomly by the normal distribution $\mathcal{N}(0,1)$, and the columns
of $\mathbf{A}$ are normalized to have a unit $\ell_{2}$-norm. The
density (the proportion of nonzero elements) of the sparse vector
$\mathbf{x}_{\textrm{true}}$ is $0.01$. The vector $\mathbf{b}$
is generated as $\mathbf{b}=(\mathbf{A}\mathbf{x}_{\textrm{true}})^{2}$.
The regularization parameter $\mu$ is set to $\mu=0.05\left\Vert \mathbf{A}^{T}\mathbf{b}\right\Vert _{\infty}$,
which allows $\mathbf{x}_{\textrm{true}}$ to be recovered to a high
accuracy.

We compare the following two instances of the proposed inexact BSCA
framework (cf. Algorithm \ref{alg:block-Successive-approximation-method-inexact}):
\begin{itemize}
\item BSCA algorithm with partial linearization approximation (Algorithm
\ref{alg:Quadratic-Inverse}), referred to as ``BSCA''. Several
variants are considered, with different number of block variables
$(K=1,2,10)$ and inner-layer iterations $(\overline{\tau}_{t}=1,10)$;
\item BSCA algorithm with quadratic approximation (cf. (\ref{eq:approximation-quadratic})),
referred to as ``BGD'' (block gradient descent). Several variants
with different number of block variables $(K=2,10)$ are considered.
The approximation subproblem has a closed-form solution and thus an
additional inner layer is not needed.
\end{itemize}
The simulation results in terms of the achieved objective value versus
the number of (outer-layer) iterations and the CPU time are shown
in Fig. \ref{fig:QI}(a)(c) and Fig. \ref{fig:QI}(b)(d), respectively.
Note that the iterations in Fig. \ref{fig:QI}(a)(c) are normalized
by the number of blocks, that is, in one iteration, all block variables
are updated once by the cyclic update rule. All algorithms start with
the same random initial point, and the stepsize is determined by the
exact line search. The quadratic regularization gain is $c_{k}^{t}=10^{-4}$
in both BSCA and BGD.

In Fig. \ref{fig:QI}(a)-(b), we investigate the impact of the number
of blocks $K$, whereas all (inexact) BSCA algorithms have the same
inner-layer iterations. We choose 10 inner-layer iterations so that
the (outer-layer) approximation subproblems can be solved with a high
accuracy. Some observations are in order.

We see that all algorithms converge to the same objective value. Note
that the BSCA algorithm with $K=1$ is in fact a fully parallel SCA
algorithm (see Sec. \ref{sec:Review-of-SCA}) and thus regarded as
the benchmark algorithm.

All BSCA algorithms with different number of blocks $(K=1,2,10)$
exhibit similar performance, in terms of both the number of iterations
and the CPU time. In practice, the number of blocks can be determined
adaptively based on the problem size and memory/computational capability
of the existing hardware. Therefore, the BSCA algorithm can solve
a much larger problem than the standard fully parallel SCA algorithm
does. In contrast, the effect of the number of blocks is more notable
for BGD algorithms.

A comparison of BSCA algorithms and BGD algorithms in Fig. \ref{fig:QI}(a)
reveals that BSCA algorithms need much fewer iterations to converge.
This consolidates the intuition that exploiting more problem structure
in the partial linearization approximate leads to faster convergence
than the general-purpose quadratic approximation (cf. the discussion
after (\ref{eq:approximation-partial-linearization})). As we see
from Fig. \ref{fig:QI}(b), this is however at the expense of more
CPU time, as the iteration complexity increases.

In Fig. \ref{fig:QI}(c)-(d) we investigate the impact of different
inner-layer iterations. Some observations are in order.

On the one hand, Fig. \ref{fig:QI}(c) shows that the BSCA with 2
blocks converges in fewer iterations when the number of inner-layer
iterations is $\overline{\tau}_{k}^{t}=10$ than when $\overline{\tau}_{k}^{t}=1$.
On the other hand, it is not surprising to see from Fig. \ref{fig:QI}(d)
that more inner-layer iterations increase the overall CPU time.

When the number of blocks is 10, the BSCA with 1 inner-layer iteration
converges in about the same number of iterations as 10 inner-layer
iterations, but its CPU time is much smaller. Hence it is not always
necessary to solve the (outer-layer) approximation subproblems with
a high accuracy.

The BSCA algorithm with a single inner-layer iteration converges in
fewer iterations than their BGD counterpart, illustrating again the
effectiveness of the partial linearization approximation that exploits
the problem structure. Furthermore, the BSCA with 10 blocks and 1
inner-layer iteration converges in about the same CPU time as the
BGD, making the inexact BSCA desirable in both the number of iterations
and the CPU time.

The BSCA with 2 blocks and 10 inner-layer iterations converges in
roughly the same number of iterations as BSCA with 10 blocks (and
either 1 or 10 inner-layer iterations). Note that at each iteration,
all blocks are updated once in the cyclic order, and after each block
update, the value of previous blocks should be passed to the next
block. This implies that the BSCA with 2 blocks and 10 inner-layer
iterations requires a smaller communication frequency than BSCA with
10 blocks.

From these observations we can see that there is no single winner.
The most suitable algorithm depends on the application and design
objective (for example, CPU time, the number of parallel processors,
the inter-communication), and it would be beneficial to incorporate
the application-specific knowledge into the algorithmic design. The
proposed algorithm is flexible enough to address different tradeoffs.

\begin{figure}[t]
\center

\includegraphics[scale=0.65]{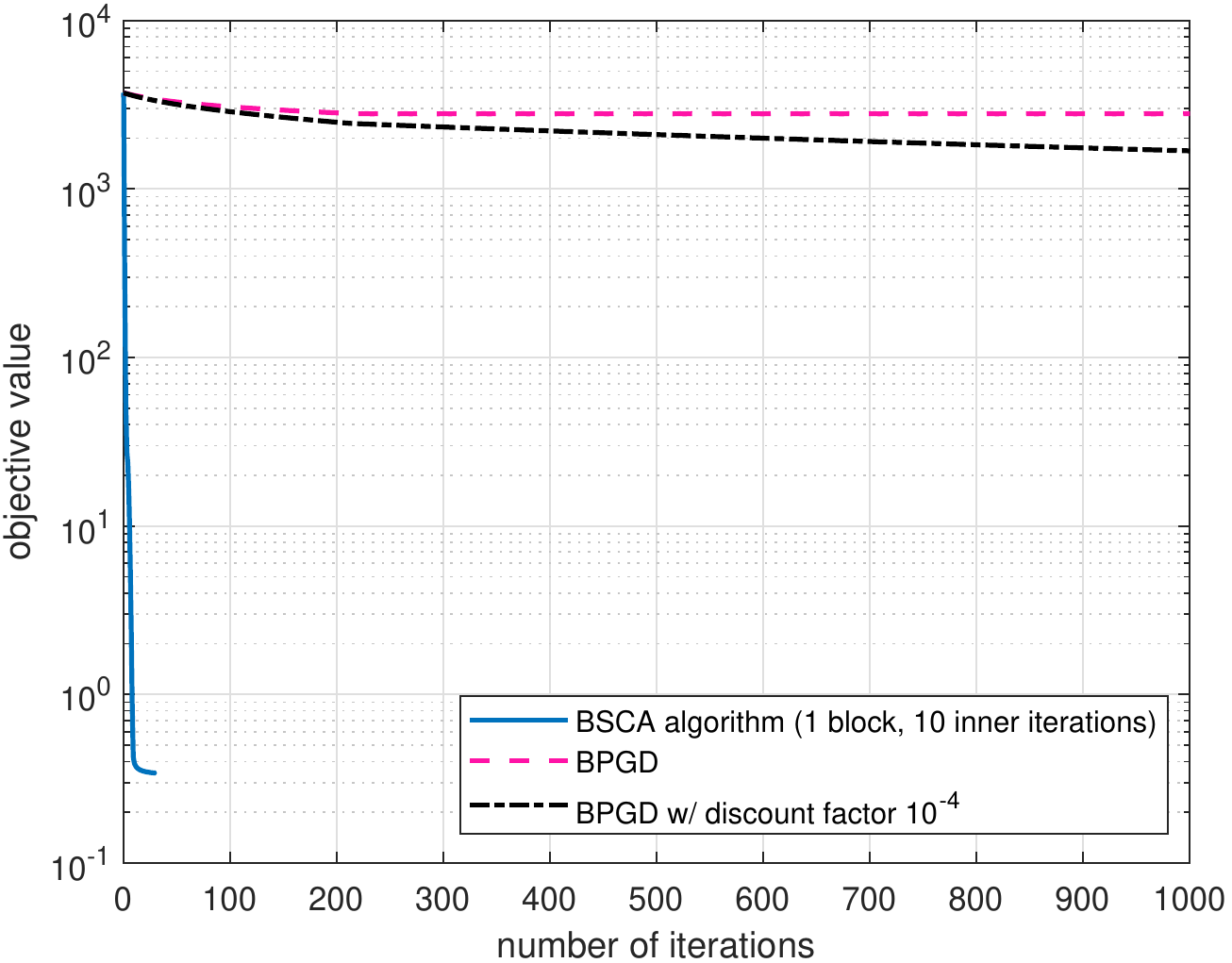}\caption{\label{fig:Bregman-PGD}Phase retrieval: BSCA and BPGD in terms of
the number of iterations}
\end{figure}

We also compare the proposed algorithm with the Bergman proximal gradient
descent (BPGD) algorithm proposed in \cite{Bolte2018}. The BPGD algorithm
extends the classical descent lemma by using non-Euclidean distances
of Bregman type. The central step is to find a constant $L$ and a
convex distance function $h$ such that both $Lh+f$ and $Lh-f$ are
convex. Particularly for the phase retrieval problem (\ref{eq:QI-formulation}),
the Bregman-based proximal gradient step at each iteration consists
of minimizing a global upper bound of the objective function $f+g$,
\begin{equation}
\min_{\mathbf{x}}\left(\frac{1}{L}\nabla f(\mathbf{x}^{t})-\nabla h(\mathbf{x}^{t})\right)^{T}\mathbf{x}+h(\mathbf{x})+\frac{1}{L}g(\mathbf{x}),\label{eq:Bregman-proximal-1}
\end{equation}
where $h(\mathbf{x})=\frac{1}{4}\left\Vert \mathbf{x}\right\Vert ^{4}+\frac{1}{2}\left\Vert \mathbf{x}\right\Vert ^{2}$
and
\begin{equation}
L=\sum_{n=1}^{N}\left(3\left\Vert \mathbf{a}_{n}\right\Vert ^{4}+\left\Vert \mathbf{a}_{n}\right\Vert ^{2}y_{n}\right).\label{eq:Bregman-proximal-2}
\end{equation}
We can see from (\ref{eq:Bregman-proximal-1}) that the value of $L$
is essential in the convergence speed: a larger $L$ indicates a less
dominating role of the function of interest $f+g$ (compared with
the distance function $h$), and thus slower convergence. The theoretical
bound (\ref{eq:Bregman-proximal-2}) usually tends to be overly conservative,
and we see from Fig. \ref{fig:Bregman-PGD} that the BPGD algorithm
converges in many more iterations than the proposed BSCA algorithm\footnote{The complexity per iteration of the BSCA and the BPGD algorithms are
comparable: both involve a soft-thresholding operator and finding
the zero of a three-order polynomial.}. Furthermore, as shown in Fig. \ref{fig:Bregman-PGD}, even if the
theoretical bound in (\ref{eq:Bregman-proximal-2}) is discounted
by a factor of $10^{-4}$, that is, $L=10^{-4}\sum_{n=1}^{N}\left(3\left\Vert \mathbf{a}_{n}\right\Vert ^{4}+\left\Vert \mathbf{a}_{n}\right\Vert ^{2}y_{n}\right)$,
convergence of BPGD (with discount factor $10^{-4}$) is still observed
in the numerical tests. Finding an appropriate value of $L$ that
yields fast convergence is a difficult task on its own. In contrast,
the proposed BSCA algorithm does not have any hyperparameters and
it thus leads to robust performance for different problem setups.
We finally note that the BPGD algorithm does not allow block updates.

\section{\label{sec:Concluding-Remarks}Concluding Remarks}

In this paper, we proposed a block successive convex approximation
algorithm for nonsmooth nonconvex optimization problems. The proposed
algorithm partitions the whole set of variables into blocks which
are updated sequentially and the dimension of each block can be adopted
to the hardware at hand. At each iteration, a block variable is selected
and updated by solving an approximation subproblem with respect to
that block variable. Compared with state-of-the-art algorithms, the
proposed algorithm has several attractive features, namely, i) high
flexibility, as the approximation function only needs to be strictly
convex and it does not have to be a global upper bound of the original
function; ii) fast convergence, as the approximation function can
be tailored to the problem at hand and the stepsize is calculated
by the line search; iii) low complexity, as the approximation subproblems
usually admit a closed-form solution and the line search scheme is
carried out over a properly constructed differentiable function; iv)
guaranteed convergence of a subsequence to a stationary point, even
when the approximation subproblem is solved inexactly and the objective
function does not have a Lipschitz continuous gradient. These attractive
features are illustrated by two applications in network anomaly detection
and phase retrieval, both theoretically and numerically.\textcolor{white}{\cite{bertsekas-ndp}
\cite{Rockafellar1998}}

\appendices{}

\section{\label{sec:Proof-of-Theorem-BSCA}Proof of Theorem \ref{thm:block-SCA-convergence}}
\begin{IEEEproof}
[Proof for the cyclic update]As the exact line search yields a larger
decrease in the objective function value than the successive line
search at each iteration, we prove the theorem without loss of generality
(w.l.o.g.) for the case where the stepsizes are calculated by the
successive line search.

Consider a limit point $\mathbf{x}^{\star}$ of the sequence $\{\mathbf{x}^{t}\}$
and a subsequence $\{\mathbf{x}^{t}\}_{t\in\mathcal{T}}$ converging
to $\mathbf{x}^{\star}$. Since $\{f(\mathbf{x}^{t})+g(\mathbf{x}^{t})\}$
is a monotonically decreasing sequence which is bounded from below, \newcounter{MYtempeqncnt} \begin{figure*}[t] \normalsize \setcounter{MYtempeqncnt}{\value{equation}} \setcounter{equation}{42} \vspace*{4pt}
\[
\frac{f(\mathbf{x}_{1}^{t}+\beta^{m_{t}-1}\left\Vert \triangle\mathbf{x}_{1}^{t}\right\Vert \frac{(\mathbb{B}_{1}\mathbf{x}^{t}-\mathbf{x}_{1}^{t})}{\left\Vert \triangle\mathbf{x}_{1}^{t}\right\Vert },\mathbf{x}_{-1}^{t})-f(\mathbf{x}_{1}^{t},\mathbf{x}_{-1}^{t})+\beta^{m_{t}-1}\left\Vert \triangle\mathbf{x}_{1}^{t}\right\Vert \frac{(g_{1}(\mathbb{B}_{1}\mathbf{x}^{t})-g_{1}(\mathbf{x}_{1}^{t}))}{\left\Vert \triangle\mathbf{x}_{1}^{t}\right\Vert }}{\beta^{m_{t}-1}\left\Vert \triangle\mathbf{x}_{1}^{t}\right\Vert }>\alpha\frac{d_{1}(\mathbf{x}^{t})}{\left\Vert \triangle\mathbf{x}_{1}^{t}\right\Vert }.
\]
\setcounter{equation}{\value{MYtempeqncnt}} \hrulefill  \end{figure*}
\begin{align*}
\lim_{\mathcal{T}\ni t\rightarrow\infty}f(\mathbf{x}^{t})+g(\mathbf{x}^{t}) & =\lim_{\mathcal{T}\ni t\rightarrow\infty}f(\mathbf{x}^{t+1})+g(\mathbf{x}^{t+1})\\
 & =f(\mathbf{x}^{\star})+g(\mathbf{x}^{\star}).
\end{align*}
By further restricting to a subsequence if necessary, we can assume
w.l.o.g. that in the subsequence $\{\mathbf{x}^{t}\}_{t\in\mathcal{T}}$
the first block is updated. It follows from the definition of the
successive line search that for all $t\in\mathcal{T}$:
\[
f(\mathbf{x}^{t+1})+g(\mathbf{x}^{t+1})-(f(\mathbf{x}^{t})+g(\mathbf{x}^{t}))\leq\alpha\beta^{m_{t}}d_{1}(\mathbf{x}^{t})\leq0,
\]
and thus
\begin{equation}
\lim_{\mathcal{T}\ni t\rightarrow\infty}\beta^{m_{t}}d_{1}(\mathbf{x}^{t})=0.\label{eq:convergence-1}
\end{equation}

From (\ref{eq:convergence-1}) we claim that
\begin{equation}
\lim_{\mathcal{T}\ni t\rightarrow\infty}\mathbb{B}_{1}\mathbf{x}^{t}-\mathbf{x}_{1}^{t}=\mathbf{0}.\label{eq:convergence-1.1}
\end{equation}
To show this, we first assume the contrary: there exists a $\delta\in(0,1)$
and a $\bar{t}$ such that
\begin{equation}
\left\Vert \mathbb{B}_{1}\mathbf{x}^{t}-\mathbf{x}_{1}^{t}\right\Vert \geq\delta,\forall\mathcal{T}\ni t\geq\bar{t}.\label{eq:convergence-lower-bound}
\end{equation}
Then (\ref{eq:convergence-1}) can be rewritten as
\begin{equation}
\lim_{\mathcal{T}\ni t\rightarrow\infty}\beta^{m_{t}}\left\Vert \triangle\mathbf{x}_{1}^{t}\right\Vert \frac{d_{1}(\mathbf{x}^{t})}{\left\Vert \triangle\mathbf{x}_{1}^{t}\right\Vert }=0,\label{eq:convergence-2}
\end{equation}
where
\[
\triangle\mathbf{x}_{1}^{t}\triangleq\left[\begin{array}{c}
\mathbb{B}_{1}\mathbf{x}^{t}\\
g_{1}(\mathbb{B}_{1}\mathbf{x}^{t})
\end{array}\right]-\left[\begin{array}{c}
\mathbf{x}_{1}^{t}\\
g_{1}(\mathbf{x}_{1}^{t})
\end{array}\right].
\]

Define
\[
\mathbf{z}_{1}^{t}\triangleq\frac{\triangle\mathbf{x}_{1}^{t}}{\left\Vert \triangle\mathbf{x}_{1}^{t}\right\Vert }.
\]
Since $\left\Vert \mathbf{z}_{1}^{t}\right\Vert =1$, by further restricting
to a subsequence if necessary, we assume the limit point of the sequence
$\{\mathbf{z}_{1}^{t}\}_{t\in\mathcal{T}}$ is $\mathbf{z}_{1}^{\star}=(\mathbf{z}_{x_{1}}^{\star},z_{g_{1}}^{\star})$
such that
\begin{align*}
\mathbf{z}_{x_{1}}^{\star} & =\lim_{\mathcal{T}\ni t\rightarrow\infty}\frac{\mathbb{B}_{1}\mathbf{x}^{t}-\mathbf{x}_{1}^{t}}{\left\Vert \triangle\mathbf{x}^{t}\right\Vert },\\
z_{g_{1}}^{\star} & =\lim_{\mathcal{T}\ni t\rightarrow\infty}\frac{g_{1}(\mathbb{B}_{1}\mathbf{x}^{t})-g_{1}(\mathbf{x}_{1}^{t})}{\left\Vert \triangle\mathbf{x}^{t}\right\Vert }.
\end{align*}
As $d_{1}(\mathbf{x}^{t})=(\mathbb{B}_{1}\mathbf{x}^{t}-\mathbf{x}_{1}^{t})^{T}\nabla_{1}f(\mathbf{x}^{t})+g_{1}(\mathbb{B}_{1}\mathbf{x}^{t})-g_{1}(\mathbf{x}_{1}^{t})$
and $\nabla f$ and $g$ are continuous functions,
\[
\lim_{\mathcal{T}\ni t\rightarrow\infty}\frac{d_{1}(\mathbf{x}^{t})}{\left\Vert \triangle\mathbf{x}_{1}^{t}\right\Vert }=\nabla_{1}f(\mathbf{x}^{\star})^{T}\mathbf{z}_{x_{1}}^{\star}+z_{g_{1}}^{\star}.
\]

There are two cases implied by (\ref{eq:convergence-2}) and we show
that neither of them is true.

\textbf{Case A: }The first case implied by (\ref{eq:convergence-2})
is that $\lim_{\mathcal{T}\ni t\rightarrow\infty}d_{1}(\mathbf{x}^{t})/\left\Vert \triangle\mathbf{x}_{1}^{t}\right\Vert =0$,
that is,
\begin{equation}
\lim_{\mathcal{T}\ni t\rightarrow\infty}\frac{d_{1}(\mathbf{x}^{t})}{\left\Vert \triangle\mathbf{x}_{1}^{t}\right\Vert }=\nabla_{1}f(\mathbf{x}^{\star})^{T}\mathbf{z}_{x_{1}}^{\star}+z_{g_{1}}^{\star}=0.\label{eq:convergence-case-1.1}
\end{equation}
Note that $\mathbf{z}_{x_{1}}^{\star}\neq\mathbf{0}$; otherwise it
implies $z_{g_{1}}^{\star}=0$, and this would contradict the fact
that $\left\Vert \mathbf{z}_{1}^{\star}\right\Vert =1$.

Since $\delta/\left\Vert \triangle\mathbf{x}_{1}^{t}\right\Vert \leq1$
and $\mathcal{X}_{1}$ is a closed and convex set, the limit point
of the following sequence is contained in $\mathcal{X}_{1}$:
\[
\lim_{\mathcal{T}\ni t\rightarrow\infty}\mathbf{x}_{1}^{t}+\frac{\delta}{\left\Vert \triangle\mathbf{x}_{1}^{t}\right\Vert }(\mathbb{B}_{1}\mathbf{x}^{t}-\mathbf{x}_{1}^{t})=\mathbf{x}_{1}^{\star}+\delta\mathbf{z}_{x_{1}}^{\star}\in\mathcal{X}.
\]
Applying the strict convexity of $\widetilde{f}(\mathbf{x}_{1};\mathbf{y})$
in $\mathbf{x}_{1}$ for any given $\mathbf{y}$, we readily obtain
\begin{align}
 & \widetilde{f}\left(\mathbf{x}_{1}^{\star}+\delta\mathbf{z}_{x_{1}}^{\star};\mathbf{x}^{\star}\right)+g_{1}(\mathbf{x}_{1}^{\star})+\delta z_{g_{1}}^{\star}\nonumber \\
>\; & \widetilde{f}(\mathbf{x}_{1}^{\star};\mathbf{x}^{\star})+\delta\nabla\widetilde{f}_{1}(\mathbf{x}_{1}^{\star};\mathbf{x}^{\star})^{T}\mathbf{z}_{x_{1}}^{\star}+g_{1}(\mathbf{x}_{1}^{\star})+\delta z_{g_{1}}^{\star}\nonumber \\
=\; & \widetilde{f}(\mathbf{x}_{1}^{\star};\mathbf{x}^{\star})+g_{1}(\mathbf{x}_{1}^{\star})+\delta(\nabla_{1}f(\mathbf{x}_{1}^{\star})^{T}\mathbf{z}_{x_{1}}^{\star}+z_{g_{1}}^{\star})\nonumber \\
=\; & \widetilde{f}(\mathbf{x}_{1}^{\star};\mathbf{x}^{\star})+g_{1}(\mathbf{x}_{1}^{\star}),\label{eq:convergence-case1-conclusion1}
\end{align}
where the equality in (\ref{eq:convergence-case1-conclusion1}) follows
from (\ref{eq:convergence-case-1.1}).

On the other hand, due to the convexity of $\widetilde{f}(\mathbf{x}_{1};\mathbf{y})$
in $\mathbf{x}_{1}$ for any given $\mathbf{y}$, we have
\begin{align*}
 & \widetilde{f}\left(\mathbf{x}_{1}^{t}+\frac{\delta}{\left\Vert \triangle\mathbf{x}_{1}^{t}\right\Vert }(\mathbb{B}_{1}\mathbf{x}^{t}-\mathbf{x}_{1}^{t});\mathbf{x}^{t}\right)\\
 & +\frac{\delta}{\left\Vert \triangle\mathbf{x}^{t}\right\Vert }g_{1}(\mathbb{B}_{1}\mathbf{x}^{t})+\left(1-\frac{\delta}{\left\Vert \triangle\mathbf{x}^{t}\right\Vert }\right)g_{1}(\mathbf{x}_{1}^{t})\\
\leq\; & \frac{\delta}{\left\Vert \triangle\mathbf{x}_{1}^{t}\right\Vert }(\widetilde{f}(\mathbb{B}_{1}\mathbf{x}^{t};\mathbf{x}^{t})+g_{1}(\mathbb{B}_{1}\mathbf{x}^{t}))\\
 & +\left(1-\frac{\delta}{\left\Vert \triangle\mathbf{x}_{1}^{t}\right\Vert }\right)(\widetilde{f}(\mathbf{x}_{1}^{t};\mathbf{x}^{t})+g_{1}(\mathbf{x}_{1}^{t}))\\
\leq\; & \widetilde{f}(\mathbf{x}_{1}^{t};\mathbf{x}^{t})+g_{1}(\mathbf{x}_{1}^{t}),
\end{align*}
where the last inequality comes from the optimality of $\mathbb{B}_{1}\mathbf{x}^{t}$.
Taking limit of the above inequality we obtain
\begin{equation}
\widetilde{f}(\mathbf{x}_{1}^{\star}+\delta\mathbf{z}_{1}^{\star};\mathbf{x}^{\star})+g_{1}(\mathbf{x}_{1}^{\star})+\delta z_{g_{1}}^{\star}\leq\widetilde{f}(\mathbf{x}_{1}^{\star};\mathbf{x}^{t})+g_{1}(\mathbf{x}_{1}^{\star}),\label{eq:convergence-case1-conclusion2}
\end{equation}
which contradicts (\ref{eq:convergence-case1-conclusion1}). Therefore
(\ref{eq:convergence-case-1.1}) cannot be true and
\begin{equation}
\lim_{\mathcal{T}\ni t\rightarrow\infty}\frac{d_{1}(\mathbf{x}^{t})}{\left\Vert \triangle\mathbf{x}_{1}^{t}\right\Vert }<0.\label{eq:convergence-case-1.2}
\end{equation}

\textbf{Case B: }The second case implied by (\ref{eq:convergence-2})
(and (\ref{eq:convergence-case-1.2})) is that
\begin{equation}
\underset{\mathcal{T}\ni t\rightarrow\infty}{\lim\sup}\;\beta^{m_{t}}\left\Vert \triangle\mathbf{x}_{1}^{t}\right\Vert =0.\label{eq:convergence-case2.0}
\end{equation}
Since $\beta^{m_{t}}\left\Vert \triangle\mathbf{x}_{1}^{t}\right\Vert \geq0$,
(\ref{eq:convergence-case2.0}) is equivalent to
\begin{equation}
\lim_{\mathcal{T}\ni t\rightarrow\infty}\beta^{m_{t}}\left\Vert \triangle\mathbf{x}_{1}^{t}\right\Vert =0.\label{eq:convergence-case2}
\end{equation}
This together with (\ref{eq:convergence-lower-bound}) implies that
$\beta^{m_{t}}\rightarrow0$, which further implies that there exists
$\bar{t}'$ such that for $\mathcal{T}\ni t\ge\bar{t}'$:
\[
\begin{array}{l}
{\displaystyle f(\mathbf{x}_{1}^{t}+\beta^{m_{t}-1}(\mathbb{B}_{1}\mathbf{x}^{t}-\mathbf{x}_{1}^{t}))+\beta^{m_{t}-1}(g_{1}(\mathbb{B}_{1}\mathbf{x}^{t})-g_{1}(\mathbf{x}_{1}^{t}))}\\
{\displaystyle \quad>f(\mathbf{x}_{1}^{t})+\alpha\beta^{m_{t}-1}d_{1}(\mathbf{x}^{t})}.
\end{array}
\]
Rearranging the terms we obtain the inequality at the top of the next
page. Letting $\mathcal{T}\ni t\rightarrow\infty$, we obtain
\[
\nabla_{1}f(\mathbf{x}^{\star})^{T}\mathbf{z}_{x_{1}}^{\star}+z_{g_{1}}^{\star}\geq\alpha(\nabla_{1}f(\mathbf{x}^{\star})^{T}\mathbf{z}_{x_{1}}^{\star}+z_{g_{1}}^{\star}),
\]
and thus
\begin{equation}
\nabla_{1}f(\mathbf{x}^{\star})^{T}\mathbf{z}_{x_{1}}^{\star}+z_{g_{1}}^{\star}\geq0.\label{eq:convergence-case2.1}
\end{equation}
Repeating the above steps (\ref{eq:convergence-case1-conclusion1})-(\ref{eq:convergence-case1-conclusion2})
in Case A (whereas the ``='' in (\ref{eq:convergence-case1-conclusion1})
should be replaced by ``$\geq$'' in view of (\ref{eq:convergence-case2.1}))
leads to a contradiction. Therefore (\ref{eq:convergence-1.1}) must
hold.

Now we show $\mathbf{x}^{\star}$ is a stationary point of (\ref{eq:problem-formulation}).
On the one hand, it follows from (\ref{eq:convergence-1.1}) that
\begin{equation}
\lim_{\mathcal{T}\ni t\rightarrow\infty}\mathbb{B}_{1}\mathbf{x}^{t}=\lim_{\mathcal{T}\ni t\rightarrow\infty}(\mathbb{B}_{1}\mathbf{x}^{t}-\mathbf{x}_{1}^{t}+\mathbf{x}_{1}^{t})=\mathbf{x}_{1}^{\star},\label{eq:convergence-3.0}
\end{equation}
and $\{\mathbb{B}_{1}\mathbf{x}^{t}\}_{t\in\mathcal{T}}$ is thus
bounded. On the other hand, it follows from the definition of $\mathbb{B}_{1}\mathbf{x}^{t}$
that
\[
\widetilde{f}(\mathbb{B}_{1}\mathbf{x}^{t};\mathbf{x}^{t})+g_{1}(\mathbb{B}_{1}\mathbf{x}^{t})\leq\widetilde{f}(\mathbf{x}_{1};\mathbf{x}^{t})+g_{1}(\mathbf{x}_{1}),\forall\mathbf{x}\in\mathcal{X},
\]
and thus
\begin{align*}
\widetilde{f}(\mathbf{x}_{1}^{\star};\mathbf{x}^{\star})+g_{1}(\mathbf{x}_{1}^{\star}) & =\lim_{\mathcal{T}\ni t\rightarrow\infty}\widetilde{f}(\mathbb{B}_{1}\mathbf{x}^{t};\mathbf{x}^{t})+g_{1}(\mathbb{B}_{1}\mathbf{x}^{t})\\
 & \leq\lim_{\mathcal{T}\ni t\rightarrow\infty}\widetilde{f}(\mathbf{x}_{1};\mathbf{x}^{t})+g_{1}(\mathbf{x}_{1})\\
 & =\widetilde{f}(\mathbf{x}_{1};\mathbf{x}^{\star})+g_{1}(\mathbf{x}_{1}),\forall\mathbf{x}\in\mathcal{X}.
\end{align*}
That is, $\mathbf{x}_{1}^{\star}$ is the optimal point of $\min_{\mathbf{x}_{1}\in\mathcal{X}_{1}}\widetilde{f}(\mathbf{x}_{1};\mathbf{x}^{\star})+g_{1}(\mathbf{x}_{1})$
and it satisfies the first order optimality condition: $g_{1}$ has
a subgradient $\boldsymbol{\xi}_{1}(\mathbf{x}_{1}^{\star})$ such
that
\begin{align*}
0 & \leq(\mathbf{x}_{1}-\mathbf{x}_{1}^{\star})^{T}(\nabla_{1}\widetilde{f}(\mathbf{x}_{1}^{\star};\mathbf{x}^{\star})+\boldsymbol{\xi}_{1}(\mathbf{x}_{1}^{\star}))\\
 & =(\mathbf{x}_{1}-\mathbf{x}_{1}^{\star})^{T}(\nabla_{1}f(\mathbf{x}^{\star})+\boldsymbol{\xi}_{1}(\mathbf{x}_{1}^{\star})),
\end{align*}
where the equality comes from Assumption (A3).

Furthermore, since $\lim_{\mathcal{T}\ni t\rightarrow\infty}\mathbb{B}_{1}\mathbf{x}^{t}-\mathbf{x}_{1}^{t}=\mathbf{0}$,
$\lim_{\mathcal{T}\ni t\rightarrow\infty}\mathbf{x}^{t+1}=\mathbf{x}^{\star}$.
In this subsequence $\{\mathbf{x}^{t+1}\}_{t\in\mathcal{T}}$, the
second block variable is updated and following the above line of analysis,
we can conclude that
\[
(\mathbf{x}_{2}-\mathbf{x}_{2}^{\star})^{T}(\nabla_{2}f(\mathbf{x}^{\star})+\boldsymbol{\xi}_{2}(\mathbf{x}_{2}^{\star}))\geq0,\forall\mathbf{x}_{2}\in\mathcal{X}_{2}.
\]
Repeating this process for the other block variables, we obtain for
$k=1,\ldots,K$ that
\[
(\mathbf{x}_{k}-\mathbf{x}_{k}^{\star})^{T}(\nabla_{k}f(\mathbf{x}^{\star})+\boldsymbol{\xi}_{k}(\mathbf{x}_{k}^{\star}))\geq0,\forall\mathbf{x}_{k}\in\mathcal{X}_{k}.
\]
Adding them up over $k=1,\ldots,K$, we readily see that $\mathbf{x}^{\star}$
satisfies the first order optimality condition, namely,
\[
(\mathbf{x}-\mathbf{x}^{\star})^{T}(\nabla f(\mathbf{x}^{\star})+\boldsymbol{\xi}(\mathbf{x}^{\star}))\geq0,\forall\mathbf{x}\in\mathcal{X},
\]
where $\boldsymbol{\xi}(\mathbf{x}^{\star})=(\boldsymbol{\xi}_{k}(\mathbf{x}_{k}^{\star}))_{k=1}^{K}$.
The proof is thus completed.
\end{IEEEproof}
\begin{IEEEproof}
[Proof for the random update]Suppose $\mathcal{U}^{t}$ is the set
of block variables that are updated at iteration $t$. It follows
from the update rule that
\[
f(\mathbf{x}^{t+1})+g(\mathbf{x}^{t+1})\leq f(\mathbf{x}^{t})+g(\mathbf{x}^{t})+\alpha\sum_{k\in\mathcal{U}^{t}}\beta^{m_{t}}d_{k}(\mathbf{x}^{t}).
\]
Introducing a Bernoulli random variable $(R_{j}^{t})_{j=1}^{K}$ where
$R_{j}^{t}$ is 1 if $\mathbf{x}_{j}$ is updated or 0 otherwise,
we can rewrite the above equation as
\[
f(\mathbf{x}^{t+1})+g(\mathbf{x}^{t+1})\leq f(\mathbf{x}^{t})+g(\mathbf{x}^{t})+\alpha\sum_{k=1}^{K}R_{k}^{t}\beta^{m_{t}}d_{k}(\mathbf{x}^{t}).
\]
Since $\mathbb{E}\left[R_{k}^{t}|\mathbf{x}^{t}\right]=p_{k}^{r}\geq p_{\min}$,
taking the expectation w.r.t. $(R_{j}^{t})_{j=1}^{K}$ conditioned
on $\mathbf{x}^{t}$ yields
\[
\mathbb{E}\left[f(\mathbf{x}^{t+1})+g(\mathbf{x}^{t+1})|\mathbf{x}^{t}\right]\leq f(\mathbf{x}^{t})+g(\mathbf{x}^{t})+\alpha p_{\min}\sum_{k=1}^{K}\beta^{m_{t}}d_{k}(\mathbf{x}^{t}).
\]
Thus $\{f(\mathbf{x}^{t})+g(\mathbf{x}^{t})\}$ is a supermartingale
w.r.t. the natural history. It follows from the supermartingale convergence
theorem \cite[Prop. 4.2]{bertsekas-ndp} that, with probability 1,
$\{f(\mathbf{x}^{t})+g(\mathbf{x}^{t})\}$ converges and
\begin{equation}
\sum_{t=0}^{\infty}\sum_{k=1}^{K}\beta^{m_{t}}d_{k}(\mathbf{x}^{t})>-\infty.\label{eq:convergence-4}
\end{equation}

Consider a limit point of the sequence $\{\mathbf{x}^{t}\}$ and a
subsequence $\{\mathbf{x}^{t}\}_{t\in\mathcal{T}}$ converging to
that limit point. By further restricting to a subsequence if necessary,
we assume w.l.o.g. that in the subsequence $\{\mathbf{x}^{t}\}_{t\in\mathcal{T}}$
the first block is updated. It follows from (\ref{eq:convergence-4})
that
\[
\sum_{t\in\mathcal{T}}\beta^{m_{t}}d_{1}(\mathbf{x}^{t})\geq\sum_{t=0}^{\infty}\sum_{k=1}^{K}\beta_{k}^{m_{t}}d_{k}(\mathbf{x}^{t})>-\infty,
\]
and thus
\begin{equation}
\underset{\mathcal{T}\ni t\rightarrow\infty}{\lim}\beta^{m_{t}}d_{1}(\mathbf{x}^{t})=0.\label{eq:convergence-5}
\end{equation}

We claim that (\ref{eq:convergence-5}) implies that $\lim_{\mathcal{T}\ni t\rightarrow\infty}\mathbb{B}_{1}\mathbf{x}^{t}-\mathbf{x}_{1}^{t}=\mathbf{0}$.
Similar to the proof for the cyclic update, we show this by contradiction:
we assume there exists a $\bar{t}$ such that $\left\Vert \mathbb{B}_{1}\mathbf{x}^{t}-\mathbf{x}_{1}^{t}\right\Vert \geq\delta$
for all $t\in\mathcal{T}$ and $t\geq\bar{t}$. Since the approximation
function is strictly convex, it follows from the previous steps that
{[}cf. Case A in (\ref{eq:convergence-case-1.2}){]}
\[
\lim_{\mathcal{T}\ni t\rightarrow\infty}\frac{d_{1}(\mathbf{x}^{t})}{\left\Vert \triangle\mathbf{x}_{1}^{t}\right\Vert }<0
\]
and $\lim\sup_{\mathcal{T}\ni t\rightarrow\infty}\beta^{m_{t}}\left\Vert \triangle\mathbf{x}_{1}^{t}\right\Vert =0$
{[}cf. Case B in (\ref{eq:convergence-case2}){]}. This further implies
that there exists a subsequence $\{\beta^{m_{t}}\}_{t\in\mathcal{T}_{s}}$
with $\mathcal{T}_{s}\subseteq\mathcal{T}$ such that $\lim_{\mathcal{T}_{s}\ni t\rightarrow\infty}\beta^{m_{t}}=0$.
This statement, however, cannot be true (see Case B of the previous
proof). Therefore $\lim_{\mathcal{T}\ni t\rightarrow\infty}\mathbb{B}_{1}\mathbf{x}^{t}-\mathbf{x}_{1}^{t}=\mathbf{0}$.

To conclude the proof, we need to show that the limit point of $\{\mathbf{x}_{t}\}$
is a stationary point. This can be proved by following the same line
of analysis in the previous proof ((\ref{eq:convergence-3.0}) and
onwards). The proof is thus completed.
\end{IEEEproof}

\section{\label{sec:Proof-of-Theorem-inexact-BSCA}Proof of Theorem \ref{thm:inexact-SCA}}
\begin{IEEEproof}
We prove the theorem w.l.o.g. for the case that only one iteration
is executed, that is, $\bar{\tau}_{t}=1$ and $\tau=0$, while the
stepsizes are calculated by the successive line search.

From (\ref{eq:descent-direction-1}) we see that the approximation
subproblem (\ref{eq:hybrid-approximate-problem}) does not have to
be solved exactly to obtain a descent direction. As a matter of fact,
repeating the same steps in (\ref{eq:descent-direction-1})-(\ref{eq:descent-direction}),
we see that $\widetilde{\mathbf{x}}_{k}^{t}-\mathbf{x}_{k}^{t}$ for
any point $\widetilde{\mathbf{x}}_{k}^{t}$ would be a descent direction
of $h(\mathbf{x}_{k},\mathbf{x}_{-k}^{t})$ at $\mathbf{x}^{t}$ as
long as
\begin{equation}
\widetilde{f}(\widetilde{\mathbf{x}}_{k}^{t};\mathbf{x}^{t})+g_{k}(\widetilde{\mathbf{x}}_{k}^{t})-(\widetilde{f}(\mathbf{x}_{k}^{t};\mathbf{x}^{t})+g_{k}(\mathbf{x}_{k}^{t}))<0.\label{eq:inexact-descent-direction}
\end{equation}

Given Assumptions (B1)-(B3), it follows from the same line of reasoning
in (\ref{eq:review-descent-direction}) (with the following notation
mapping: $\mathbf{x}^{t}\rightarrow\mathbf{x}_{k}^{t,\tau}$, $\mathbb{B}\mathbf{x}^{t}\rightarrow\mathbb{B}_{k}\overline{\mathbf{x}}^{t,\tau}$,
$\widetilde{f}(\mathbf{x};\mathbf{x}^{t})\rightarrow\widetilde{f}^{i}(\mathbf{x}_{k};\mathbf{x}_{k}^{t,\tau},\mathbf{x}^{t})$,
$f(\mathbf{x})\rightarrow\widetilde{f}(\mathbf{x}_{k};\mathbf{x}^{t})$,
$g(\mathbf{x})\rightarrow g_{k}(\mathbf{x}_{k})$) that
\begin{align}
d_{k}(\overline{\mathbf{x}}^{t,0})=\; & (\mathbb{B}_{k}\overline{\mathbf{x}}^{t,0}-\mathbf{x}_{k}^{t,0})^{T}\nabla\widetilde{f}(\mathbf{x}_{k}^{t,0};\mathbf{x}_{k}^{t})\nonumber \\
 & +g_{k}(\mathbb{B}_{k}\overline{\mathbf{x}}^{t,0})-g_{k}(\mathbf{x}_{k}^{t,0})<0.\label{eq:inexact-block-SCA-proof-0}
\end{align}
Since $\mathbf{x}_{k}^{t,1}$ is obtained by performing the successive
line search along $\widetilde{f}(\mathbf{x}_{k};\mathbf{x}^{t})+g_{k}(\mathbf{x}_{k})$
(cf. Step 2.3 of Alg. \ref{alg:block-Successive-approximation-method-inexact}),
we have
\begin{equation}
\widetilde{f}(\mathbf{x}_{k}^{t,1};\mathbf{x}^{t})+g_{k}(\mathbf{x}_{k}^{t,1})\negthickspace\leq\negthickspace\widetilde{f}(\mathbf{x}_{k}^{t,0};\mathbf{x}^{t})+g_{k}(\mathbf{x}_{k}^{t,0})+\alpha\beta^{m_{t,0}}d_{k}(\overline{\mathbf{x}}^{t,0}).\label{eq:inexact-block-SCA-proof-1}
\end{equation}

Combining (\ref{eq:inexact-block-SCA-proof-1}) with (\ref{eq:inexact-block-SCA-proof-0})
and recall $\widetilde{\mathbf{x}}_{k}^{t}=\mathbf{x}_{k}^{t,1}$
and $\mathbf{x}_{k}^{t,0}=\mathbf{x}_{k}^{t}$, we readily obtain
(\ref{eq:inexact-descent-direction}) and thus
\begin{align*}
0 & >(\widetilde{\mathbf{x}}_{k}^{t}-\mathbf{x}_{k}^{t})^{T}\nabla\widetilde{f}(\mathbf{x}_{k}^{t};\mathbf{x}^{t})+g_{k}(\widetilde{\mathbf{x}}_{k}^{t})-g_{k}(\mathbf{x}_{k}^{t})\\
 & =(\widetilde{\mathbf{x}}_{k}^{t}-\mathbf{x}_{k}^{t})^{T}\nabla_{k}f(\mathbf{x}^{t})+g_{k}(\widetilde{\mathbf{x}}_{k}^{t})-g_{k}(\mathbf{x}_{k}^{t}),
\end{align*}
where the inequality is due to the convexity of $\widetilde{f}(\mathbf{x};\mathbf{x}^{t})$
(Assumption (A1)) and the equality is due to Assumption (A3). Since
$\mathbf{x}_{k}^{t+1}$ is obtained by performing the line search
over $f+g$ along the coordinate of $\mathbf{x}_{k}$, cf. Step S3
of Alg. \ref{alg:block-Successive-approximation-method-inexact},
we have
\begin{align*}
f & (\mathbf{x}^{t+1})+g(\mathbf{x}^{t+1})-(f(\mathbf{x}^{t})+g(\mathbf{x}^{t}))\\
 & \leq\alpha\beta^{m_{t}}((\widetilde{\mathbf{x}}_{k}^{t}-\mathbf{x}_{k}^{t})^{T}\nabla_{k}f(\mathbf{x}^{t})+g_{k}(\widetilde{\mathbf{x}}_{k}^{t})-g_{k}(\mathbf{x}_{k}^{t}))\leq0,
\end{align*}

Consider a limit point $\mathbf{x}^{\star}$ and a subsequence $\{\mathbf{x}^{t}\}_{t\in\mathcal{T}}$
converging to $\mathbf{x}^{\star}$, and assume w.l.o.g. that the
first block is updated in this subsequence. Following the same line
of analysis from (\ref{eq:convergence-lower-bound}) to (\ref{eq:convergence-case2.1})
(with the notation mapping $\mathbb{B}_{k}\mathbf{x}^{t}\rightarrow\widetilde{\mathbf{x}}_{k}^{t}$),
we have
\[
\lim_{\mathcal{T}\ni t\rightarrow\infty}\widetilde{\mathbf{x}}_{1}^{t}-\mathbf{x}_{1}^{t}=\mathbf{0}.
\]
and furthermore
\begin{align}
\lim_{\mathcal{T}\ni t\rightarrow\infty}\widetilde{\mathbf{x}}_{1}^{t} & =\lim_{\mathcal{T}\ni t\rightarrow\infty}(\widetilde{\mathbf{x}}_{1}^{t}-\mathbf{x}_{1}^{t}+\mathbf{x}_{1}^{t})\nonumber \\
 & =\lim_{\mathcal{T}\ni t\rightarrow\infty}(\widetilde{\mathbf{x}}_{1}^{t}-\mathbf{x}_{1}^{t})+\lim_{\mathcal{T}\ni t\rightarrow\infty}\mathbf{x}_{1}^{t}=\mathbf{x}_{1}^{\star}.\label{eq:inexact-block-SCA-proof-4}
\end{align}

We claim that $\mathbf{x}_{1}^{\star}$ is the optimal point of the
following outer-layer approximation subproblem (cf. (\ref{eq:hybrid-approximate-problem})):
\begin{equation}
\underset{\mathbf{x}_{1}\in\mathcal{X}_{1}}{\textrm{minimize}}\;\widetilde{f}(\mathbf{x}_{1};\mathbf{x}^{\star}),\label{eq:inexact-block-SCA-proof-3}
\end{equation}
and we show this by contradiction.

First of all, $\mathbf{x}_{1}^{\star}$ is the optimal point of problem
(\ref{eq:inexact-block-SCA-proof-3}) if and only if it is an optimal
point of the following inner-layer approximation subproblem
\begin{equation}
\underset{\mathbf{x}_{1}\in\mathcal{X}_{1}}{\textrm{minimize}}\;\widetilde{f}^{i}(\mathbf{x}_{1};\mathbf{x}_{1}^{\star},\mathbf{x}^{\star}).\label{eq:inexact-block-SCA-proof-3.1}
\end{equation}
Define $\overline{\mathbf{x}}^{\star}\triangleq(\mathbf{x}_{1}^{\star},\mathbf{x}^{\star})$.
If $\mathbf{x}_{1}^{\star}$ is not the optimal point of (\ref{eq:inexact-block-SCA-proof-3.1}),
we denote as $\mathbb{B}_{1}\overline{\mathbf{x}}^{\star}$ an optimal
point of (\ref{eq:inexact-block-SCA-proof-3.1}). Then $\mathbb{B}_{1}\overline{\mathbf{x}}^{\star}\neq\mathbf{x}_{1}^{\star}$
and $\mathbb{B}_{1}\overline{\mathbf{x}}^{\star}-\mathbf{x}_{1}^{\star}$
is a descent direction of $\widetilde{f}^{i}(\mathbf{x}_{1};\mathbf{x}_{1}^{\star},\mathbf{x}^{\star})$
at $\mathbf{x}_{1}=\mathbf{x}_{1}^{\star}$, in the sense that $d_{1}(\overline{\mathbf{x}}^{\star})<0$.

Recall $\mathbf{x}_{1}^{t,0}=\mathbf{x}_{1}^{t}$ and the optimality
of $\mathbb{B}_{1}\overline{\mathbf{x}}_{1}^{t,0}$ with $\overline{\mathbf{x}}_{1}^{t,0}=(\mathbf{x}_{1}^{t,0},\mathbf{x}^{t})$,
we note that for any $\mathbf{x}_{1}\in\mathcal{X}_{1}$,
\[
\widetilde{f}^{i}(\mathbf{x}_{1};\mathbf{x}_{1}^{t},\mathbf{x}^{t})+g_{1}(\mathbf{x}_{1})\geq\widetilde{f}^{i}(\mathbb{B}_{1}\overline{\mathbf{x}}_{1}^{t,0};\mathbf{x}_{1}^{t},\mathbf{x}^{t})+g_{1}(\mathbb{B}_{1}\overline{\mathbf{x}}_{1}^{t,0}).
\]
Since $\{\mathbb{B}_{1}\overline{\mathbf{x}}_{1}^{t,0}\}_{t\in\mathcal{T}}$
is bounded by Assumption (B5), it has a convergent subsequence and
we denote its limit point as $\mathbf{y}_{1}$. Restricting to that
sequence if necessary, we have
\begin{align*}
\widetilde{f}^{i}(\mathbf{x}_{1};\mathbf{x}_{1}^{\star},\mathbf{x}^{\star})+ & g_{1}(\mathbf{x}_{1})=\lim_{\mathcal{T}\ni t\rightarrow\infty}\widetilde{f}^{i}(\mathbf{x}_{1};\mathbf{x}_{1}^{t},\mathbf{x}^{t})+g_{1}(\mathbf{x}_{1})\\
 & \geq\lim_{\mathcal{T}\ni t\rightarrow\infty}\widetilde{f}^{i}(\mathbb{B}_{1}\overline{\mathbf{x}}_{1}^{t,0};\mathbf{x}_{1}^{t},\mathbf{x}^{t})+g_{1}(\mathbb{B}_{1}\overline{\mathbf{x}}_{1}^{t,0})\\
 & =\widetilde{f}^{i}\left(\mathbf{y}_{1};\mathbf{x}_{1}^{\star},\mathbf{x}^{\star}\right)+g_{1}\left(\mathbf{y}_{1}\right),\forall\mathbf{x}_{1}\in\mathcal{X}_{1}.
\end{align*}
Therefore,
\[
\lim_{\mathcal{T}\ni t\rightarrow\infty}\left.\mathbb{B}_{1}\overline{\mathbf{x}}^{t,0}\right|_{\overline{\mathbf{x}}^{t,0}=(\mathbf{x}_{1}^{t},\mathbf{x}^{t})}=\mathbf{y}_{1}=\mathbb{B}_{1}\overline{\mathbf{x}}^{\star}\bigr|_{\overline{\mathbf{x}}^{\star}=(\mathbf{x}_{1}^{\star},\mathbf{x}^{\star})}.
\]
Since the iterative algorithm in the inner layer is executed for one
iteration only,
\begin{equation}
\widetilde{\mathbf{x}}_{1}^{t}=\mathbf{x}_{1}^{t,1}=\mathbf{x}_{1}^{t}+\gamma^{t,0}(\mathbb{B}_{1}\overline{\mathbf{x}}^{t,0}-\mathbf{x}_{1}^{t}),\label{eq:inexact-block-SCA-proof-5-1}
\end{equation}
where $\gamma^{t,0}$ is the stepsize obtained by applying successive
line search to $\widetilde{f}_{1}(\mathbf{x}_{1};\mathbf{x}^{t})+g_{1}(\mathbf{x}_{1})$.

This successive line search consists of two conceptual steps. The
first conceptual step is to identify the set of $\gamma$ such that
\[
\left\{ \gamma(\mathbf{x}^{t})\geq0\left|\begin{array}{l}
\widetilde{f}(\mathbf{x}_{1}^{t}+\gamma(\mathbf{x}^{t})(\mathbb{B}_{1}\overline{\mathbf{x}}^{t,0}-\mathbf{x}_{1}^{t});\mathbf{x}^{t})\smallskip\\
+g_{1}(\mathbf{x}_{1}^{t})+\gamma(\mathbf{x}^{t})(g_{1}(\mathbb{B}_{1}\overline{\mathbf{x}}^{t,0})-g_{1}(\mathbf{x}_{1}^{t}))\smallskip\\
=\widetilde{f}(\mathbf{x}_{1}^{t};\mathbf{x}^{t})+g_{1}(\mathbf{x}_{1}^{t})+\alpha\gamma(\mathbf{x}^{t})d(\overline{\mathbf{x}}^{t,0})
\end{array}\right.\right\} .
\]
This set is a singleton since $d(\overline{\mathbf{x}}^{t,0})<0$
and $\widetilde{f}(\mathbf{x}_{1};\mathbf{x}^{t})$ is strictly convex.
The second conceptual step is to identify the smallest nonnegative
integer $m_{t,0}$ such that $\beta^{m_{t,0}}\leq\gamma(\mathbf{x}^{t})$.
We assume w.l.o.g. that $\{\gamma(\mathbf{x}^{t})\}_{t\in\mathcal{T}}$
is bounded; otherwise $\beta^{m_{t,0}}=1$ and $\lim_{\mathcal{T}\ni t\rightarrow\infty}\beta^{m_{t,0}}=1$.
Restricting to a convergent subsequence of $\{\gamma(\mathbf{x}^{t})\}_{t\in\mathcal{T}}$
if necessary, it follows from \cite[5.8 Example]{Rockafellar1998}
that $\lim_{\mathcal{T}\ni t\rightarrow\infty}\gamma(\mathbf{x}^{t})=\gamma(\mathbf{x}^{\star})$,
where $\gamma(\mathbf{x}^{\star})>0$ satisfies
\begin{align*}
 & \widetilde{f}(\mathbf{x}_{1}^{\star}+\gamma(\mathbf{x}^{\star})(\mathbb{B}_{1}\overline{\mathbf{x}}^{\star}-\mathbf{x}_{1}^{\star});\mathbf{x}^{\star})+g_{1}(\mathbf{x}_{1}^{\star})\\
 & +\gamma(\mathbf{x}^{\star})(g_{1}(\mathbb{B}_{1}\overline{\mathbf{x}}^{\star})-g_{1}(\mathbf{x}_{1}^{\star}))\\
=\; & \widetilde{f}(\mathbf{x}_{1}^{\star};\mathbf{x}^{\star})+g_{1}(\mathbf{x}_{1}^{\star})+\alpha\gamma(\mathbf{x}^{\star})d(\overline{\mathbf{x}}^{\star}).
\end{align*}
Therefore, $\lim_{\mathcal{T}\ni t\rightarrow\infty}\beta^{m_{t,0}}=\beta^{m^{\star}}>0$,
where $m^{\star}$ is the smallest nonnegative integer such that $\beta^{m^{\star}}\leq\gamma(\mathbf{x}^{\star})$.
Taking the limit of (\ref{eq:inexact-block-SCA-proof-5-1}), we have
\begin{align}
\lim_{\mathcal{T}\ni t\rightarrow\infty}\widetilde{\mathbf{x}}_{1}^{t} & =\lim_{\mathcal{T}\ni t\rightarrow\infty}(\mathbf{x}_{1}^{t}+\beta^{m_{t,0}}(\mathbb{B}_{1}\overline{\mathbf{x}}^{t,0}-\mathbf{x}_{1}^{t}))\nonumber \\
 & =\mathbf{x}_{1}^{\star}+\beta^{m^{\star}}(\mathbb{B}_{1}\overline{\mathbf{x}}^{\star}-\mathbf{x}_{1}^{\star}).\label{eq:inexact-block-SCA-proof-5}
\end{align}
Since $\beta^{m^{\star}}>0$ in (\ref{eq:inexact-block-SCA-proof-5}),
it is also valid for the case that $\{\gamma(\mathbf{x}^{t})\}_{t\in\mathcal{T}_{s}}$
is unbounded and $\lim_{\mathcal{T}\ni t\rightarrow\infty}\beta^{m_{t,0}}=1$.

A comparison between the two equations (\ref{eq:inexact-block-SCA-proof-4})
and (\ref{eq:inexact-block-SCA-proof-5}) implies that $\mathbf{x}_{1}^{\star}=\mathbb{B}_{1}\overline{\mathbf{x}}^{\star}$,
and hence a contradiction is derived. Therefore, $\mathbf{x}_{1}^{\star}$
is the optimal point of (\ref{eq:inexact-block-SCA-proof-3}). By
following the same line of analysis of the proof of Theorem \ref{thm:block-SCA-convergence}
in Appendix \ref{sec:Proof-of-Theorem-BSCA}, we can repeat the same
steps for $\mathbf{x}_{2},\mathbf{x}_{3},\ldots,\mathbf{x}_{K}$.
Therefore $\mathbf{x}^{\star}$ is a stationary point of (\ref{eq:problem-formulation})
and the proof is completed.
\end{IEEEproof}

\end{document}